\documentclass[12pt]{article}
\usepackage{amsmath,amsthm,amssymb,mathtools}
\allowdisplaybreaks
\usepackage{enumitem}
\usepackage{etoolbox}
\usepackage[colorlinks,
linkcolor=red,
anchorcolor=blue,
citecolor=magenta
]{hyperref}
\hypersetup{pdftitle={Uniform Exponential Stability for Finite Energy Weak Solutions},pdfauthor={Xin Song, Chengfeng Xiong},pdfsubject={Barotropic compressible Navier--Stokes equations with large potential forces}}

\theoremstyle{plain}
\newtheorem{theorem}{Theorem}[section]
\newtheorem{lemma}[theorem]{Lemma}
\newtheorem{corollary}[theorem]{Corollary}
\newtheorem{definition}{Definition}[section]
\theoremstyle{definition}
\newtheorem{remark}{Remark}[section]
\theoremstyle{plain}

\numberwithin{equation}{section}

\newcommand{\R}{\mathbb R}
\newcommand{\dd}{\,\mathrm d}

\newcommand{\supp}{\operatorname{supp}}

\newcommand{\Data}{\mathfrak D}

\title{Uniform Exponential Stability
for Finite Energy Weak Solutions to the Barotropic
Compressible Navier--Stokes Equations with
Large Potential External Forces}
\author{%
\normalsize Xin Song\textsuperscript{a},
Chengfeng Xiong\textsuperscript{b,c}\thanks{Corresponding author.\newline
\protect\raggedright Email addresses: \mbox{songxin2357@163.com} (X. Song),
\mbox{xiongcf4math@126.com} (C.F. Xiong).}\\
\small a. School of Mathematics and Computer Sciences,\\
\small Nanchang University, Nanchang 330031, P. R. China\\
\small b. Institute of Applied Mathematics,\\
\small Academy of Mathematics and Systems Science (AMSS),\\
\small Chinese Academy of Sciences, Beijing 100190, P. R. China\\
\small c. China Mobile Internet Company Limited,\\
\small Guangzhou 510640, P. R. China}
\date{}

\makeatletter
\patchcmd{\@title}{\\}{ }{}{}
\patchcmd{\@title}{\\}{ }{}{}
\patchcmd{\@title}{\\}{ }{}{}
\patchcmd{\@maketitle}{\large}{\normalsize}{}{}
\let\XiongOriginalMaketitle\maketitle
\renewcommand{\maketitle}{\setcounter{footnote}{1}\XiongOriginalMaketitle\vspace{-.7cm}}
\renewcommand{\@biblabel}[1]{[#1]}
\makeatother

\newif\ifXiongKeywordStyle
\AfterEndEnvironment{abstract}{%
  \par\begingroup\XiongKeywordStyletrue\footnotesize\let\noindent\indent
}
\pretocmd{\section}{\ifXiongKeywordStyle\par\endgroup\fi}{}{}

\BeforeBeginEnvironment{thebibliography}{%
  \begingroup\footnotesize
  \let\textsc\textnormal
  \let\emph\textnormal
  \let\textbf\textnormal
}
\AfterEndEnvironment{thebibliography}{\endgroup}

\begin{document}
\maketitle
\begin{abstract}
We study the large-time behavior of Lions--Feireisl finite energy weak
solutions to the barotropic compressible Navier--Stokes equations in
three-dimensional bounded domains with Dirichlet boundary conditions and
large potential external forces. For $\gamma>3/2$, we prove that these weak solutions converge exponentially
to the steady state as time tends to infinity, provided that the steady
density is bounded away from zero. The main novelty of our analysis is a compactness argument based on
Young measures. Furthermore, we prove that the Young measures generated
by the density are Dirac measures. Together with the extra integrability of the density, this yields strong
convergence to the nonconstant steady density and uniform smallness of
the energy $w$ for sufficiently large time. We use Taylor expansion and the Bogovskii
operator to construct a Lyapunov functional, and use this smallness of
$w$ to obtain exponential decay with a time $T_0$, which is
independent of the particular weak solution, thereby
answering the question raised by Kang, Peng and Xiong
in \cite[Remark~1.1]{KPX}.
Corresponding decay estimates for the density, momentum, kinetic energy
and the viscous dissipation integrated from large times to infinity also
hold, allowing vacuum.
\end{abstract}


\section{Introduction and main results}\label{eqb-setting}
We study the three-dimensional barotropic compressible Navier--Stokes system
\begin{equation}\label{eqb-system}
\begin{cases}
\partial_t\rho+\operatorname{div}(\rho u)=0,\\
\partial_t(\rho u)+\operatorname{div}(\rho u\otimes u)+\nabla P(\rho)
 =\mu\Delta u+(\lambda+\mu)\nabla\operatorname{div}u+\rho\nabla F,\\
u|_{\partial\Omega}=0.
\end{cases}
\end{equation}
Here $\rho\ge0$ and $u\in\R^3$ denote the density and velocity, respectively. Furthermore, $\Omega\subset\R^3$ is a bounded connected
$C^{2+\kappa}$ domain for some $\kappa>0$. The constant viscosity coefficients $\mu$ and $\lambda$ satisfy
\begin{equation}\label{eqb-parameter-conditions}
\mu>0,\qquad 2\mu+3\lambda\ge0.
\end{equation}
The potential $F\in W^{1,\infty}(\Omega)$ is independent of time and the pressure is
\begin{equation}\label{eqb-pressure-law}
P(\rho)=\rho^\gamma,\qquad \gamma>3/2.
\end{equation}
We further impose the initial conditions
\begin{equation}\label{eqb-initial-conditions}
\rho(x,0)=\rho_0,\qquad (\rho u)(x,0)=m_0.
\end{equation}

We use the continuous representative of $F$ and write
\begin{equation}\label{eqb-potential-bounds}
F_*:=\inf_\Omega F\quad\text{and}\quad F^*:=\sup_\Omega F.
\end{equation}
and for $f\in L^1(\Omega)$, we denote
\[
\langle f\rangle:=\frac1{|\Omega|}\int_\Omega f\,dx
\]
for its spatial average over $\Omega$.

The well-posedness and large-time behavior of compressible Navier--Stokes
equations have been studied extensively. Now, we list some related important results in the multi-dimensional case. We first recall some results
on global weak solutions in several space dimensions. For the pressure
$P(\rho)=\rho^\gamma$, Lions \cite{Lions1998} established global existence
of weak solutions in three  dimensional space where $\gamma\ge9/5$.
Feireisl--Novotn\'y--Petzeltov\'a \cite{FNP2001} extended this result to
$\gamma>3/2$, which includes the case $\gamma=5/3$ for a monoatomic gas.
For spherically symmetric initial data, Jiang--Zhang \cite{JZ2001}
proved global existence of weak solutions to the Cauchy problem for
$\gamma>1$.

In discussing compactness of finite energy weak solutions, Feireisl
\cite{Feireisl2001} explained that, for $\gamma\ge9/5$, square
integrability of the densities and their weak limit implies that the
limiting pair satisfies the continuity equation in the renormalized
sense of DiPerna--Lions \cite{DiPernaLions1989}. He noted that this
conclusion is not immediate under the weaker restriction $\gamma>3/2$.
His approach uses density cut-offs to establish the renormalized
continuity equation for the limit and then to prove strong convergence
of the density.

For initial data close to a constant state with positive density,
Matsumura--Nishida \cite{MN1980} obtained global classical solutions.
Hoff \cite{Hoff1995} treated discontinuous initial data under a
smallness assumption. Huang--Li--Xin \cite{HLX2012} established global
existence and uniqueness of classical solutions to the three-dimensional
isentropic Cauchy problem for smooth initial data with small energy,
allowing large oscillations and vacuum. In bounded domains, Cai--Li
\cite{CL2023} obtained global weak and classical solutions with
Navier-type slip boundary conditions for suitably small initial energy,
allowing vacuum.
For the
Cauchy problem with large potential forces, Li--Zhang--Zhao
\cite{LZZ2015} proved global existence and uniqueness of strong
solutions with small initial energy, allowing large oscillations and
interior vacuum, provided that the steady density is bounded away from
zero.

For the stationary problem with potential forces, Feireisl--Petzeltov\'a \cite{FP1998}
established uniqueness of the stationary density with prescribed mass
when the upper level sets of the external potential are connected.
Erban \cite{Erban2001} examined the
multiplicity of stationary solutions with prescribed mass and energy.

Concerning the large-time behavior, Matsumura--Padula \cite{MP1992} studied the stability of stationary
solutions in bounded domains under large potential forces.
For weak solutions, Feireisl--Petzeltov\'a
\cite{FP1999} studied stabilization of compressible flows with
Dirichlet boundary conditions. Novotn\'y--Stra\v skraba \cite{NS2000}
investigated stabilization under potential forces for boundary
conditions that admit a Poincar\'e inequality, including Dirichlet
conditions. Feireisl--Petzeltov\'a analyzed the propagation of density
oscillations in \cite{FP2001compact} and established the existence of
bounded absorbing sets for the isentropic system when $\gamma>5/3$ in
\cite{FP2001absorbing}.

For exponential decay of barotropic flows, Padula \cite{Padula1999} proved exponential
stability of the rest state in a bounded domain for a class of weak
perturbations with density bounded away from zero. For the isentropic
equations without external forces, Fang--Zi--Zhang \cite{FZZ2012} used
the Bogovskii operator to construct a Lyapunov functional and prove
exponential convergence in $L^2$, assuming an upper bound on the
density but allowing vacuum. Peng--Shi--Wu \cite{PengShi} established
exponential decay for Lions--Feireisl finite energy weak solutions in
three-dimensional bounded domains without additional upper or lower
bounds on the density. Their argument uses the extra integrability of
the density and a suitable Lyapunov functional.
Xu--Zhong \cite{XZ2026} studied exponential
convergence in $L^2$ in three-dimensional bounded domains with slip
boundary conditions and large potential forces, using small initial
potential energy and a uniform upper bound on the density.

Young measures were introduced by L.C. Young to treat variational
problems without classical minimizers, as recalled by Ball
\cite{Ball1989}. Ball gave a version of the fundamental theorem for
Young measures that represents weak limits of nonlinear functions,
under suitable hypotheses, by a family of probability measures.
Gallenm\"uller--Wagner--Wiedemann \cite{GWW2023} reviewed their use
in studying the limits of oscillating sequences and measure-valued
formulations of fluid equations.

Kang, Peng and Xiong \cite{KPX} proved exponential decay to the
equilibrium state for Lions--Feireisl finite energy weak solutions
in three-dimensional bounded domains with large potential external
forces. However, the time $T_0$ in their decay estimate depends on the
particular weak solution. Whether $T_0$ can be chosen independently
of the weak solution was left open in \cite[Remark~1.1]{KPX}. In this
paper, we give an affirmative answer to this question under the
assumptions of Theorem~\ref{thm:main}. We use Young measures to obtain compactness of density sequence . An identity
for the pressure and truncated densities shows that the Young
measures generated by the density are Dirac measures. Together with
the extra integrability of the density, this gives strong convergence
to the steady density. This compactness result yields a common time,
depending only on the prescribed data, after which the relative energy $w$,
defined in \eqref{eqb-energy}, is small. We then use this smallness
and a Lyapunov functional constructed with the Bogovskii operator
to prove exponential decay with $T_0$ and the decay constants
depending only on the prescribed data.

Let $m=\int_\Omega\rho_0\,dx$ and assume
\begin{equation}\label{eqb-mass-condition}
m>S(F_*),
\end{equation}
where
\[
S(k):=\int_\Omega\left[\frac{\gamma-1}{\gamma}(F(x)-k)_+\right]^{1/(\gamma-1)}dx,
\]
and $(a)_+=\max\{a,0\}$. The unique nonnegative equilibrium density is given by
\begin{equation}\label{eqb-equilibrium}
\rho_s(x)
 =\left[\frac{\gamma-1}{\gamma}(F(x)-k_0)\right]^{1/(\gamma-1)}.
\end{equation}
where $k_0<F_*$ is uniquely determined by $S(k_0)=m$.
In particular, $\rho_s\in W^{1,\infty}(\Omega)$ and $0<\underline{\rho_s}\le \rho_s\le \overline{\rho_s}$, where
\begin{equation}\label{eqb-r-bounds}
\begin{aligned}
\underline{\rho_s}&:=\left[\frac{\gamma-1}{\gamma}(F_*-k_0)\right]^{1/(\gamma-1)}>0,\\
\overline{\rho_s}&:=\left[\frac{\gamma-1}{\gamma}(F^*-k_0)\right]^{1/(\gamma-1)}<\infty.
\end{aligned}
\end{equation}

We now introduce the definition of finite energy weak solutions, following
Feireisl--Novotn\'y--Petzeltov\'a \cite{FNP2001}.

\begin{definition}\label{def:weak-solution}
A pair $(\rho,u)$ is a finite energy weak solution of
\eqref{eqb-system} on $(0,\infty)\times\Omega$ if the following conditions hold:
\begin{itemize}
\item The density is nonnegative and
\[
\rho\in L^\infty(0,\infty;L^\gamma(\Omega)),\qquad
u\in L^2(0,\infty;H^1_0(\Omega;\R^3)).
\]
\item Both equations in \eqref{eqb-system} hold in
$\mathcal D'((0,\infty)\times\Omega)$. The continuity equation also holds
in $\mathcal D'((0,\infty)\times\R^3)$ after extending $\rho,u$ by zero
outside $\Omega$.
\item The mechanical energy
\begin{equation}\label{eqb-mechanical-energy}
\mathcal E(t):=\int_\Omega\left(
\frac12\rho|u|^2+\frac{\rho^\gamma}{\gamma-1}-\rho F\right)dx
\end{equation}
belongs to $L^1_{\mathrm{loc}}(0,\infty)$ and satisfies
\begin{equation}\label{eqb-mechanical-distribution}
\begin{aligned}
-\int_0^\infty\varphi'(t)\mathcal E(t)\,dt
 +\int_0^\infty\varphi(t)\int_\Omega
 \bigl(\mu|\nabla u|^2+(\lambda+\mu)(\operatorname{div}u)^2\bigr)dx\,dt
\le0
\end{aligned}
\end{equation}
for every nonnegative $\varphi\in C_c^\infty(0,\infty)$.
\item The continuity equation holds in the renormalized sense:
\begin{equation}\label{eqb-renormalized}
\partial_t b(\rho)+\operatorname{div}(b(\rho)u)
 +[\rho b'(\rho)-b(\rho)]\operatorname{div}u=0
\end{equation}
in $\mathcal D'((0,\infty)\times\Omega)$ for every $b\in C^1(\R)$ such
that $b'(z)=0$ for $|z|\ge M$, where $M$ may depend on $b$.
\end{itemize}
\end{definition}

Next, we prescribe initial data satisfying the compatibility conditions in
\cite{FNP2001}:
\begin{equation}\label{eqb-initial-data}
\rho_0\in L^\gamma(\Omega),\quad \rho_0\ge0,\quad
m_0=0\ \text{a.e. on }\{\rho_0=0\},\quad
\frac{|m_0|^2}{\rho_0}\in L^1(\Omega).
\end{equation}
The quotient is equal to zero on vacuum.
For the solutions considered in this paper, the weak continuity properties
recalled in \cite{FNP2001} take the following form: for every finite $T>0$,
\begin{equation}\label{eqb-weak-continuity}
\rho\in C_w([0,T];L^\gamma(\Omega)),\quad
\rho u\in C_w([0,T];L^{2\gamma/(\gamma+1)}(\Omega;\R^3)),
\end{equation}
where $C_w$ denotes continuity for the weak topology. It follows from \eqref{eqb-mechanical-distribution} that the energy inequality
\begin{equation}\label{eqb-relative-initial-energy}
\int_\Omega\left(G(\rho,\rho_s)+\frac12\rho|u|^2\right)dx
 +\int_0^t\!\int_\Omega
 \bigl[\mu|\nabla u|^2+(\lambda+\mu)(\operatorname{div}u)^2\bigr]dx\,ds
 \le w_0
\end{equation}
holds for almost every $t\ge0$, where
\[
w_0:=\int_\Omega
 \left(\frac{|m_0|^2}{2\rho_0}+G(\rho_0,\rho_s)\right)dx
\quad\text{and}\quad
G(\rho,\rho_s):=\rho\int_{\rho_s}^{\rho}
 \frac{h^\gamma-\rho_s^\gamma}{h^2}\,dh.
\]
By Taylor's formula with integral remainder,
for $s,z>0$ we obtain
\begin{equation}\label{eqb-G}
\begin{aligned}
G(s,z)=H(s)-H(z)-H'(z)(s-z)=\frac{s^\gamma}{\gamma-1}-\frac{\gamma z^{\gamma-1}s}{\gamma-1}+z^\gamma,
\end{aligned}
\end{equation}
where $H(s):=\frac{s^\gamma}{\gamma-1}$.

We define the relative energy with respect to $(\rho_s,0)$
(see \cite{FJN2012,Sueur2014}) by
\begin{equation}\label{eqb-energy}
w(t):=\int_\Omega\left(\frac12\rho|u|^2+G(\rho,\rho_s)\right)dx,
\end{equation}
and set $\Data:=(\Omega,\gamma,\mu,\lambda,F,m,w_0)$. We are now ready to formulate our main result.
\begin{theorem}\label{thm:main}
Let $\Omega\subset\R^3$ be a bounded, connected domain
of class $C^{2+\kappa}$ for some $\kappa>0$.
Suppose that $F\in W^{1,\infty}(\Omega)$ and \eqref{eqb-parameter-conditions}--\eqref{eqb-pressure-law} hold. Moreover, the initial data $(\rho_0,m_0)$ satisfy
\eqref{eqb-initial-data} and \eqref{eqb-mass-condition}.
Then there exist constants $T_0\ge0$ and $c>0$, depending only on $\Data$,
such that every finite-energy weak solution of \eqref{eqb-system}
satisfies, for almost every $t>T_0$,
\begin{align}
w(t)&\le3w_0e^{-c(t-T_0)},\label{eq:main-decay}\\
\|\rho(t)-\rho_s\|_{L^\gamma(\Omega)}
&\le C_d(3w_0)^\alpha e^{-\alpha c(t-T_0)},\label{eq:density-rate}
\end{align}
and
\begin{equation}\label{eq:momentum-rate}
\begin{aligned}
&\|\sqrt\rho u(t)\|_{L^2(\Omega)}
+\|\rho u(t)\|_{L^{2\gamma/(\gamma+1)}(\Omega)}
\le\sqrt{6w_0}
\left[1+\left(\overline{\rho_s}|\Omega|^{1/\gamma}+C_d\right)^{1/2}\right]
e^{-\frac c2(t-T_0)},
\end{aligned}
\end{equation}
where $\alpha:=1/\max\{2,\gamma\}$ and
\[
C_d:=
\begin{cases}
\displaystyle
\left[|\Omega|^{1-\gamma/2}
\left(\frac2\gamma\right)^{\gamma/2}
(2\overline{\rho_s})^{\gamma(2-\gamma)/2}
+\frac{2^\gamma(\gamma-1)}{2^\gamma-\gamma-1}\right]^{1/\gamma},
&3/2<\gamma\le2,\\[8pt]
\displaystyle
\left[\frac8{3\gamma}
\left(\frac{4\overline{\rho_s}}{\underline{\rho_s}}\right)^{\gamma-2}
+\frac{2^\gamma(\gamma-1)}{2^\gamma-\gamma-1}\right]^{1/\gamma},
&\gamma>2.
\end{cases}
\]
Moreover, for every $t>T_0$,
\begin{equation}\label{eq:tail-dissipation}
\int_t^\infty\|\nabla u(s)\|_{L^2(\Omega)}^2\,ds
\le\frac{3w_0}{\mu}e^{-c(t-T_0)}.
\end{equation}
The constants $T_0$ and $c$ are independent of the particular weak solution
$(\rho,u)$.
\end{theorem}

\begin{remark}\label{rem:weak-solution-existence}
The global finite energy weak solutions considered in
Theorem~\ref{thm:main} exist by \cite[Theorem~1.1, Section~5]{FNP2001}.
Moreover, they satisfy
$\rho\in L^{\gamma+\theta_0}((0,T)\times\Omega)$ for every $T>0$,
where $\theta_0\in(0,(2\gamma-3)/3]$. Theorem~\ref{thm:main} shows that
these weak solutions decay exponentially to the equilibrium state.
The initial density may contain vacuum, and no smallness condition
is imposed on the initial energy.
\end{remark}

\begin{remark}\label{rem:equilibrium}
Under the mass condition \eqref{eqb-mass-condition}, the equilibrium
density given by \eqref{eqb-equilibrium} is unique and satisfies
$\rho_s\in W^{1,q}(\Omega)$ and $\inf_\Omega\rho_s>0$ for
$F\in W^{1,q}(\Omega)$, $3<q\le\infty$; see \cite[Lemma~2.5]{KPX}.
The positive lower bound on $\rho_s$ is used for the factor
$\rho_s^{-1}$ in the test function \eqref{eq:weighted-field}; see also
\cite[Remark~1.2]{KPX}.
\end{remark}

\begin{remark}\label{rem:density-exponent}
For $3/2<\gamma<2$, the exponent $\alpha=1/2$ in
\eqref{eq:density-rate} follows by combining \eqref{eq:main-decay}
with \eqref{eq:cor-density-control} once $w(t)\le1$.
In general, the estimate
\[
\|\rho-\rho_s\|_{L^\gamma(\Omega)}
\le C\left(\int_\Omega G(\rho,\rho_s)\,dx\right)^{1/\gamma}
\]
cannot hold with a constant $C$ independent of $\rho$, even for
positive densities with the same total mass $m$.
Indeed, take a nonzero $f\in L^\infty(\Omega)$ with
$\int_\Omega f\,dx=0$ and set $\rho_\varepsilon=\rho_s+\varepsilon f$.
Then $\int_\Omega\rho_\varepsilon\,dx=\int_\Omega\rho_s\,dx$,
and $\rho_\varepsilon>0$ for all sufficiently small $\varepsilon>0$,
since $\inf_\Omega\rho_s>0$. Taylor expansion gives
\[
\lim_{\varepsilon\to0^+}
\frac{\displaystyle\int_\Omega G(\rho_\varepsilon,\rho_s)\,dx}
{\varepsilon^2}
=\frac\gamma2\int_\Omega\rho_s^{\gamma-2}f^2\,dx>0.
\]
Hence, for sufficiently small $\varepsilon>0$, we have
\[
0<\int_\Omega G(\rho_\varepsilon,\rho_s)\,dx
\le\gamma\varepsilon^2\int_\Omega\rho_s^{\gamma-2}f^2\,dx.
\]
Since $\|\rho_\varepsilon-\rho_s\|_{L^\gamma(\Omega)}
=\varepsilon\|f\|_{L^\gamma(\Omega)}$ and $(2-\gamma)/\gamma>0$, we obtain
\[
\begin{aligned}
\frac{\|\rho_\varepsilon-\rho_s\|_{L^\gamma(\Omega)}}
{\left(\displaystyle\int_\Omega G(\rho_\varepsilon,\rho_s)\,dx\right)^{1/\gamma}}
&\ge
\frac{\|f\|_{L^\gamma(\Omega)}}
{\left(\gamma\displaystyle\int_\Omega\rho_s^{\gamma-2}f^2\,dx\right)^{1/\gamma}}
\frac1{\varepsilon^{(2-\gamma)/\gamma}}\\
&\longrightarrow+\infty
\qquad\text{as }\varepsilon\to0^+.
\end{aligned}
\]
Thus the estimate with power $1/\gamma$ cannot hold with a constant
independent of $\varepsilon$. 
\end{remark}

\begin{remark}\label{rem:large-time-role}
Exponential decay for Lions--Feireisl weak solutions in the absence
of external forces was proved in \cite{PengShi}. For large potential
forces and slip boundary conditions, Xu--Zhong \cite{XZ2026} obtained exponential
convergence under a uniform upper bound on the density and small initial
potential energy. Theorem~\ref{thm:main} treats the no-slip problem in
the finite energy class, with no pointwise upper bound on the
density and no smallness requirement on $w_0$.

The constants $T_0$ and $c$ in Theorem~\ref{thm:main} depend
only on $\Data$ and are independent of the particular weak solution
$(\rho,u)$. For fixed $\Data$, the estimates
\eqref{eq:main-decay}--\eqref{eq:tail-dissipation} hold with the same
$T_0$ and $c$ for every weak solution satisfying the assumptions of
Theorem~\ref{thm:main}. 
\end{remark}

We next recall the definition of Young measures used in the density
compactness argument; see \cite{GWW2023}.

\begin{definition}\label{def:young-measure}
Let $Q\subset\R^4$ be bounded and Lebesgue measurable. A Young measure
with values in $[0,\infty)$ is a family of probability measures
$(\nu_{t,x})_{(t,x)\in Q}$ on $[0,\infty)$ such that
\[
(t,x)\longmapsto\langle\nu_{t,x},f\rangle
\quad\text{is measurable for every }f\in C_0(\R),
\]
where $C_0(\R)$ denotes the continuous functions tending to zero at
infinity, and
\[
\langle\nu_{t,x},f\rangle:=\int_{[0,\infty)}f(z)\,d\nu_{t,x}(z).
\]
 A sequence of nonnegative measurable functions
$\rho_n$ on $Q$ is said to generate $(\nu_{t,x})$ if
\[
\langle\delta_{\rho_n},f\rangle\overset{*}{\rightharpoonup}\langle\nu_{t,x},f\rangle
\quad\text{in }L^\infty(Q)\quad\text{for every }f\in C_0(\R),
\]
where $\delta_a$ is the Dirac measure at $a\ge0$.
\end{definition}

\begin{remark}\label{rem:young-oscillations}
In Lemma~\ref{lem:repair-compact}, we use Young measures to identify
the weak limit of $\rho_n^\gamma$ when the viscous dissipation tends
to zero. Weak convergence of $\rho_n$ alone does not imply
$\widetilde\rho=\bar\rho^\gamma$. By
Lemma~\ref{lem:young-measure-convergence}, the weak limits satisfy
\[
\bar\rho=\langle\nu_{t,x},z\rangle
\quad\text{and}\quad
\widetilde\rho=\langle\nu_{t,x},z^\gamma\rangle.
\]
Using \eqref{eq:repair-covariance} for all the truncations $\zeta_k$,
we prove that $\nu_{t,x}=\delta_{\bar\rho(t,x)}$ almost everywhere.
Thus $\widetilde\rho=\bar\rho^\gamma$, and
Lemma~\ref{lem:young-measure-convergence} gives strong convergence
of $\rho_n$ to $\bar\rho$ in $L^\gamma(Q)$.
The limit satisfies the stationary equation and has total mass $m$.
The uniqueness in Remark~\ref{rem:equilibrium} then yields
$\bar\rho=\rho_s$.
\end{remark}

\begin{remark}\label{rem:young-stabilization}
The higher integrability estimate \eqref{eqb-hi-window} gives uniform
integrability of $\{\rho_n^\gamma\}$. Once the Young measures are shown
to be Dirac measures, Lemma~\ref{lem:young-measure-convergence} yields
strong convergence of the density in $L^\gamma(Q)$.
On unit time intervals where the viscous dissipation tends to zero,
this convergence to $\rho_s$ and the kinetic energy estimate in
Lemma~\ref{lem:repair-compact} give
$\int_0^1w_n(t)\,dt\to0$, as stated in \eqref{eq:repair-strong}.
Since $w_n$ is non-increasing, the energy at the end of each interval
also tends to zero. This proves that the lower bound $\eta_\varepsilon$
defined in \eqref{eq:uniform-cost} is positive. Combining this lower bound with the energy inequality
yields a time
$T_\varepsilon$, depending only on $\Data$ and $\varepsilon$, after
which $w(t)\le\varepsilon$ for every weak solution under consideration.
For sufficiently small $\varepsilon$, we use this smallness and a
Lyapunov functional constructed with the Bogovskii operator to prove
exponential decay, as shown in
Lemma~\ref{lem:small-energy-decay}.
\end{remark}

We now outline the proof of Theorem~\ref{thm:main}.
First, for every $\varepsilon>0$, we prove that $w(t)\le\varepsilon$
for all sufficiently large $t$, with a time depending only on $\Data$
and $\varepsilon$. We consider sequences on unit time intervals where
the viscous dissipation tends to zero. The bounds
\eqref{eqb-hi-window} and \eqref{eqb-hi-basic} are uniform for the
time translations used in the proof. Together with H\"older's,
Poincar\'e's and Sobolev's inequalities, they give convergence of the
velocities and convective terms to zero. To treat the pressure, we
apply the renormalized continuity equation to the truncated densities
and estimate the time derivatives of
$\mathcal A[\chi\zeta_k(\rho_n)]$. These estimates, the bounds for
Riesz transforms, and compact embedding allow us
to use these functions, multiplied by smooth cut-off functions, as
test functions in the momentum equation; see
Lemma~\ref{lem:local-test-compactness}.
Comparing the resulting identity with the equation satisfied by the
weak limits gives \eqref{eq:repair-covariance}. The Young measure
representation in Lemma~\ref{lem:young-measure-convergence}, based on
\cite{Ball1989}, expresses these limits as integrals against
$\nu_{t,x}$. Using \eqref{eq:repair-covariance} for all the truncations
$\zeta_k$, we prove that $\nu_{t,x}$ is a Dirac measure. The higher
integrability bound and Lemma~\ref{lem:young-measure-convergence}
then give strong convergence in $L^\gamma(Q)$. The stationary
equation, the prescribed mass, and uniqueness in
Remark~\ref{rem:equilibrium} identify the limit  $\rho_s$.
Together with convergence to zero of the kinetic energy integrated
over the unit time interval, this proves \eqref{eq:repair-strong}.
This compactness result and the monotonicity of $w$ imply that the
dissipation has a positive lower bound on each unit time interval
contained in $[1,\infty)$ whose final energy exceeds $\varepsilon$; see \eqref{eq:uniform-cost}.
This bound depends only on $\Data$ and $\varepsilon$.
Summing it over intervals on which the final energy
exceeds $\varepsilon$ and using \eqref{eqb-relative-initial-energy}
gives $T_\varepsilon$ in \eqref{eq:uniform-time}.
Thus $w(t)\le\varepsilon$ for $t\ge T_\varepsilon$, independently
of the particular weak solution, as proved in
Lemma~\ref{thm:repair-uniform}.
For exponential decay, we use the Bogovskii test function
\eqref{eq:weighted-field} and the estimate \eqref{eq:absorb-decay}
from \cite[proof of Proposition~3.1, (3.48)]{KPX}.
Choosing $\varepsilon$ sufficiently small allows the first term on
the right-hand side of \eqref{eq:absorb-decay} to be absorbed into
the left-hand side.
Combining the resulting estimate with the mechanical energy
inequality, we obtain exponential decay for the Lyapunov functional
in \eqref{eq:corrected-energy-comparison}. Its equivalence to $w$
then gives \eqref{eq:main-decay}.
The estimates for $G$, H\"older's inequality, and the integrated
energy inequality then yield
\eqref{eq:density-rate}--\eqref{eq:tail-dissipation}.

The remainder of the paper is organized as follows.
Section~\ref{eqb-coercivity-section} collects the energy and higher
integrability estimates and the Young measure convergence criterion.
The main results, Theorem~\ref{thm:main}, are proved in
Section~\ref{sec:small-energy-decay}.

\section{Preliminaries}
\label{eqb-coercivity-section}\label{eqb-hi-section}
Choose $0\le\eta\in C_c^\infty(\R^3)$ with $\int_{\R^3}\eta\,dx=1$.
For $\varepsilon>0$ and $f\in L^1_{\mathrm{loc}}(\R^3)$, set
\[
\eta_\varepsilon(x):=\varepsilon^{-3}\eta(x/\varepsilon),
\qquad [f]_\varepsilon:=\eta_\varepsilon*f.
\]
We first recall the regularized continuity equation;
see \cite[Corollary~2.3]{KPX}.

\begin{lemma}\label{lem:mollified-continuity}
Let $(\rho,u)$ be a finite energy weak solution, extended by zero outside
$\Omega$. Then, for $0<\theta\le\gamma/2$ and $I=[\alpha,\beta]\subset(0,\infty)$,
\begin{equation}\label{eq:mollified-continuity}
\partial_t[\rho^\theta]_\varepsilon
+\operatorname{div}([\rho^\theta]_\varepsilon u)
=(1-\theta)[\rho^\theta\operatorname{div}u]_\varepsilon+r_\varepsilon \quad\text{a.e. on}\ I\times\Omega,
\end{equation}
 where
$
r_\varepsilon
:=\operatorname{div}([\rho^\theta]_\varepsilon u)
-[\operatorname{div}(\rho^\theta u)]_\varepsilon.
$
Moreover,
\begin{equation}\label{eq:mollifier-commutator-limit}
r_\varepsilon\longrightarrow0
\quad\text{in }L^2(I;L^q(\Omega)),\quad \text{ for }
 1\le q\le\frac{2\gamma}{\gamma+2\theta}.
\end{equation}
\end{lemma}

The Bogovskii operator first constructed by Bogovskii
\cite{Bogovskii1980} provides the basis for our construction of test
functions for the momentum equation.

\begin{lemma}\label{lem:bogovskii}
There exists a linear operator $\mathcal B$ such that, for every
$1<p<\infty$, the following properties hold:
\begin{enumerate}[label=\textup{(\roman*)},ref=(\roman*)]
\item\label{item:B-bound}
$\mathcal B:L^p_0(\Omega)\to
W^{1,p}_0(\Omega;\R^3)$ is bounded, and
$\|\mathcal B[f]\|_{W^{1,p}(\Omega)}
\le C(\Omega,p)\|f\|_{L^p(\Omega)}$, where
$L^p_0(\Omega):=\{f\in L^p(\Omega):\int_\Omega f\,dx=0\}$.
\item\label{item:B-divergence}
For $f\in L^p_0(\Omega)$, $v=\mathcal B[f]$ satisfies
$\operatorname{div}v=f$ in $\Omega$ and $v|_{\partial\Omega}=0$.
\item\label{item:B-divergence-data}
For $1<r<\infty$, suppose $f=\operatorname{div}g$ with
$g\in L^r(\Omega;\R^3)$ and $g\cdot n|_{\partial\Omega}=0$.
Then $\|\mathcal B[f]\|_{L^r(\Omega)}
\le C(\Omega,r)\|g\|_{L^r(\Omega)}$.
\end{enumerate}
\end{lemma}

We now present the properties of $G$ for the convenience of the subsequent discussion.

\begin{lemma}\label{eqb-coercivity}
For $s,z\ge0$, $G(s,z)$ has the following properties:
\begin{enumerate}[label=\textup{(\roman*)},ref=(\roman*)]
\item\label{item:G-positive}
$G(s,z)\ge0,\ G(s,z)=0\iff s=z.$
\item\label{item:G-upper}
$G(s,z)\le\frac{s^\gamma}{\gamma-1}+z^\gamma.$
\item\label{item:G-tail}
$s\ge2z\ \Longrightarrow\ G(s,z)\ge c_\gamma s^\gamma,
\text{ where } c_\gamma:=\frac{2^\gamma-\gamma-1}{2^\gamma(\gamma-1)}>0.$
\item\label{item:G-quadratic}
For $z\in[\underline{\rho_s},\overline{\rho_s}]$ and $0\le s\le2\overline{\rho_s}$,
$G(s,z)\ge q_\gamma(s-z)^2$, where $q_\gamma=\dfrac\gamma2(2\overline{\rho_s})^{\gamma-2}$ for $\gamma\in (1,2]$ and $q_\gamma=\dfrac{3\gamma}{8}(\underline{\rho_s}/2)^{\gamma-2}$ for $\gamma\in (2,+\infty)$.
\item\label{item:G-weighted}
For fixed $0<\theta<\gamma/2$, there are constants
$C_a=C_a(\Data,\theta)>0$ and $C_A=C_A(\Data,\theta)>0$ such that,
for $s\ge0$ and $\underline{\rho_s}\le z\le\overline{\rho_s}$,
\begin{equation}\label{eq:pointwise-decay}
z^{-\theta}(s^\gamma-z^\gamma)(s^\theta-z^\theta)\ge C_a G(s,z)
\quad\text{and}\quad
|s^\theta-z^\theta|^2\le C_A G(s,z).
\end{equation}
\end{enumerate}
\end{lemma}
\begin{proof}
It is easy to check that \ref{item:G-positive} and \ref{item:G-upper} hold by \eqref{eqb-G}. To prove \ref{item:G-tail}, we first assume $s>0$. For $0<z\le s$,
$\partial_zG(s,z)=\gamma z^{\gamma-2}(z-s)\le0$. Thus for $0\le z\le s/2$ we obtain
\[
G(s,z)\ge G(s,s/2)
=\frac{2^\gamma-\gamma-1}{2^\gamma(\gamma-1)}s^\gamma
=c_\gamma s^\gamma.
\]
If $s=0$, the condition $s\ge2z$ follows $z=0$, which implies same inequality holds.

To prove \ref{item:G-quadratic}, let $z\in[\underline{\rho_s},\overline{\rho_s}]$ and
$0<s\le2\overline{\rho_s}$. Taylor's formula with integral remainder yields
\[
G(s,z)=\gamma(s-z)^2
\int_0^1(1-t)[z+t(s-z)]^{\gamma-2}\,dt.
\]
If $1<\gamma\le2$, then $0<z+t(s-z)\le2\overline{\rho_s}$ for
$0\le t\le1$. Hence
\[
G(s,z)\ge\gamma(2\overline{\rho_s})^{\gamma-2}(s-z)^2
\int_0^1(1-t)\,dt
=\frac\gamma2(2\overline{\rho_s})^{\gamma-2}(s-z)^2.
\]
If $\gamma>2$, then $z+t(s-z)\ge z/2\ge\underline{\rho_s}/2$ for
$0\le t\le1/2$ which gives
\[
G(s,z)\ge\gamma(\underline{\rho_s}/2)^{\gamma-2}(s-z)^2
\int_0^{1/2}(1-t)\,dt
=\frac{3\gamma}{8}(\underline{\rho_s}/2)^{\gamma-2}(s-z)^2.
\]
These inequalities extend to $s=0$ by continuity.

To prove \ref{item:G-weighted}, we first apply Taylor's formula to obtain
\[
\lim_{s\to z}
\frac{z^{-\theta}(s^\gamma-z^\gamma)(s^\theta-z^\theta)}{G(s,z)}
=2\theta,
\qquad
\lim_{s\to z}\frac{|s^\theta-z^\theta|^2}{G(s,z)}
=\frac{2\theta^2}{\gamma}z^{2\theta-\gamma}.
\]
At the endpoints, it follows from \eqref{eqb-G} that
\[
\lim_{s\to0^+}
\frac{z^{-\theta}(s^\gamma-z^\gamma)(s^\theta-z^\theta)}{G(s,z)}=1,
\qquad
\lim_{s\to\infty}
\frac{z^{-\theta}(s^\gamma-z^\gamma)(s^\theta-z^\theta)}{G(s,z)}
=+\infty,
\]
and
\[
\lim_{s\to0^+}\frac{|s^\theta-z^\theta|^2}{G(s,z)}=z^{2\theta-\gamma},
\qquad
\lim_{s\to\infty}\frac{|s^\theta-z^\theta|^2}{G(s,z)}=0,
\]
since $2\theta<\gamma$. Thus
\ref{item:G-weighted} holds.
\end{proof}

By Lemma~\ref{eqb-coercivity}, we obtain the following estimates
for the density.
\begin{corollary}\label{cor:energy-modulus-decay}
There is $C_d=C_d(\Data)>0$ such that
\begin{equation}\label{eq:cor-density-control}
\|\rho-\rho_s\|_{L^\gamma(\Omega)}
\le C_d\max\left\{
\left(\int_\Omega G(\rho,\rho_s)\,dx\right)^{1/2},
\left(\int_\Omega G(\rho,\rho_s)\,dx\right)^{1/\gamma}\right\}.
\end{equation}
and
\[
\|\rho^\theta-\rho_s^\theta\|_{L^{\gamma/\theta}(\Omega)}
\le\|\rho-\rho_s\|_{L^\gamma(\Omega)}^\theta,\qquad 0<\theta\le1.
\]
\end{corollary}
\begin{proof}
By Lemma~\ref{eqb-coercivity}~\ref{item:G-quadratic} and
\ref{item:G-tail}, we respectively have
\[
\begin{aligned}
\int_{\{\rho\le2\overline{\rho_s}\}}|\rho-\rho_s|^2\,dx
&\le q_\gamma^{-1}\int_\Omega G(\rho,\rho_s)\,dx,\\
\int_{\{\rho>2\overline{\rho_s}\}}|\rho-\rho_s|^\gamma\,dx
&\le\int_{\{\rho>2\overline{\rho_s}\}}\rho^\gamma\,dx
\le c_\gamma^{-1}\int_\Omega G(\rho,\rho_s)\,dx.
\end{aligned}
\]
For $3/2<\gamma<2$, it follows from H\"older's inequality that
\[
\begin{aligned}
\int_{\{\rho\le2\overline{\rho_s}\}}|\rho-\rho_s|^\gamma\,dx
&\le|\Omega|^{1-\gamma/2}
\left(\int_{\{\rho\le2\overline{\rho_s}\}}|\rho-\rho_s|^2\,dx\right)^{\gamma/2}\\
&\le|\Omega|^{1-\gamma/2}q_\gamma^{-\gamma/2}
\left(\int_\Omega G(\rho,\rho_s)\,dx\right)^{\gamma/2},
\end{aligned}
\]
which implies
\[
\begin{aligned}
\|\rho-\rho_s\|_{L^\gamma(\Omega)}^\gamma
&\le|\Omega|^{1-\gamma/2}q_\gamma^{-\gamma/2}
\left(\int_\Omega G(\rho,\rho_s)\,dx\right)^{\gamma/2}
+c_\gamma^{-1}\int_\Omega G(\rho,\rho_s)\,dx\\
&\le\left[|\Omega|^{1-\gamma/2}q_\gamma^{-\gamma/2}
+c_\gamma^{-1}\right]
\max\left\{\left(\int_\Omega G(\rho,\rho_s)\,dx\right)^{\gamma/2},
\int_\Omega G(\rho,\rho_s)\,dx\right\}.
\end{aligned}
\]
For $\gamma\ge2$, we can apply the fact that
$|\rho-\rho_s|\le2\overline{\rho_s}$ on
$\{\rho\le2\overline{\rho_s}\}$ to get
\[
\begin{aligned}
\|\rho-\rho_s\|_{L^\gamma(\Omega)}^\gamma
&\le(2\overline{\rho_s})^{\gamma-2}
\int_{\{\rho\le2\overline{\rho_s}\}}|\rho-\rho_s|^2\,dx
+\int_{\{\rho>2\overline{\rho_s}\}}|\rho-\rho_s|^\gamma\,dx\\
&\le\left[(2\overline{\rho_s})^{\gamma-2}q_\gamma^{-1}
+c_\gamma^{-1}\right]\int_\Omega G(\rho,\rho_s)\,dx.
\end{aligned}
\]
Thus the first estimate holds with
\[
C_d:=
\begin{cases}
\displaystyle
\left[|\Omega|^{1-\gamma/2}q_\gamma^{-\gamma/2}
+c_\gamma^{-1}\right]^{1/\gamma},
&3/2<\gamma<2,\\[6pt]
\displaystyle
\left[(2\overline{\rho_s})^{\gamma-2}q_\gamma^{-1}
+c_\gamma^{-1}\right]^{1/\gamma},
&\gamma\ge2.
\end{cases}
\]

Observe that, for $0<\theta\le1$,
\[
|s^\theta-z^\theta|^{\gamma/\theta}\le|s-z|^\gamma,
\qquad s,z\ge0,
\]
which implies
\[
\begin{aligned}
\|\rho^\theta-\rho_s^\theta\|_{L^{\gamma/\theta}(\Omega)}^{\gamma/\theta}
&=\int_\Omega|\rho^\theta-\rho_s^\theta|^{\gamma/\theta}\,dx\le\int_\Omega|\rho-\rho_s|^\gamma\,dx
=\|\rho-\rho_s\|_{L^\gamma(\Omega)}^\gamma.
\end{aligned}
\]
\end{proof}

We next show that the smallness of the energy $w(t)$ implies the smallness of $\|\rho^\theta-\rho_s^\theta\|_{\gamma/\theta}$ and $\|\rho-\rho_s\|_\gamma$.
\begin{lemma}\label{eqb-threshold}
Let $\alpha=\min\{1/2,1/\gamma\}$ and $C_d$ be the constant in
Corollary~\ref{cor:energy-modulus-decay}. For fixed $\theta\in(0,\gamma/2]$,
there is $C_\theta=C_\theta(\Data,\theta)>0$ such that
$0\le w(t)\le1$ implies
\begin{equation}\label{eqb-norm-control}
\|\rho^\theta-\rho_s^\theta\|_{\gamma/\theta}\le C_\theta w(t)^{\theta/\gamma},
\qquad \|\rho-\rho_s\|_\gamma\le C_d w(t)^\alpha.
\end{equation}
In particular, for any $\eta>0$, 
we have $w(t)\le\varepsilon_0$ implies
$\|\rho^\theta-\rho_s^\theta\|_{\gamma/\theta}+\|\rho-\rho_s\|_\gamma\le\eta$,
where
\begin{equation}\label{eqb-explicit-threshold}
	\varepsilon_0:=\min\left\{1,
	\left(\frac\eta{2C_\theta}\right)^{\gamma/\theta},
	\left(\frac\eta{2C_d}\right)^{1/\alpha}\right\}>0.
\end{equation}
\end{lemma}
\begin{proof}
Let $A=\{\rho\le2\rho_s\}$ and $B=\{\rho>2\rho_s\}$, and define
\[
L_\theta=\max\left\{\frac{2(2\overline{\rho_s})^\theta}{\underline{\rho_s}},
\theta\max\{(\underline{\rho_s}/2)^{\theta-1},(2\overline{\rho_s})^{\theta-1}\}\right\}.
\]
For $0\le s\le2\overline{\rho_s}$ and $\underline{\rho_s}\le z\le\overline{\rho_s}$,
if $s<\underline{\rho_s}/2$, then we have $|s-z|\ge\frac{\underline{\rho_s}}2$, which further yields that
\[
|s^\theta-z^\theta|\le(2\overline{\rho_s})^\theta
\le\frac{2(2\overline{\rho_s})^\theta}{\underline{\rho_s}}|s-z|.
\]
If $s\ge\underline{\rho_s}/2$, we can use the mean value theorem to get
\[
|s^\theta-z^\theta|
\le\theta\max\{(\underline{\rho_s}/2)^{\theta-1},(2\overline{\rho_s})^{\theta-1}\}|s-z|.
\]
Combining these estimates, we obtain
\[
|s^\theta-z^\theta|\le L_\theta|s-z|.
\]

Since $\gamma/\theta\ge2$, it follows from Lemma~\ref{eqb-coercivity}~\ref{item:G-quadratic} that
\[
|\rho^\theta-\rho_s^\theta|^{\gamma/\theta}
\le (2\overline{\rho_s})^{\gamma-2\theta}|\rho^\theta-\rho_s^\theta|^2
\le (2\overline{\rho_s})^{\gamma-2\theta}L_\theta^2q_\gamma^{-1}G(\rho,\rho_s),
\qquad\text{on }A.
\]
Similarly, we also have
\[
|\rho^\theta-\rho_s^\theta|^{\gamma/\theta}\le\rho^\gamma
\le c_\gamma^{-1}G(\rho,\rho_s),\qquad\text{on }B.
\]
The first estimate in \eqref{eqb-norm-control} holds with
\[
C_\theta:=\left((2\overline{\rho_s})^{\gamma-2\theta}L_\theta^2q_\gamma^{-1}+c_\gamma^{-1}\right)^{\theta/\gamma}.
\]

For $0\le w(t)\le1$, it follows from \eqref{eq:cor-density-control} in
Corollary~\ref{cor:energy-modulus-decay} and \eqref{eqb-energy} that
\[
\begin{aligned}
\|\rho-\rho_s\|_{L^\gamma(\Omega)}
&\le C_d\max\left\{
\left(\int_\Omega G(\rho,\rho_s)\,dx\right)^{1/2},
\left(\int_\Omega G(\rho,\rho_s)\,dx\right)^{1/\gamma}\right\}\\
&\le C_d\max\{w(t)^{1/2},w(t)^{1/\gamma}\}
=C_dw(t)^\alpha.
\end{aligned}
\]
This proves \eqref{eqb-norm-control}.
Finally, by \eqref{eqb-explicit-threshold}, we have
\[
w(t)\le\varepsilon_0
\quad\Longrightarrow\quad
C_\theta w(t)^{\theta/\gamma}+C_d w(t)^\alpha
\le\frac{\eta}{2}+\frac{\eta}{2}=\eta,
\]
which proves the small-energy assertion.
\end{proof}

We next establish an  higher estimate of integrability  on finite time intervals.

\begin{lemma}\label{eqb-hi}
Under the hypotheses of Theorem~\ref{thm:main}, 
then there is $C$, depending only on $\Data$, such that for $1\le a<b<\infty$,
\begin{equation}\label{eqb-hi-window}
\int_a^b\int_\Omega\rho^{\gamma+\theta}\,dxdt\le C(1+b-a).
\end{equation}
where $\theta:=\min\left\{\frac{2\gamma-3}{6},\frac16\right\}$
\end{lemma}
\begin{proof}
By \eqref{eqb-relative-initial-energy} and
Lemma~\ref{eqb-coercivity}~\ref{item:G-tail},
\begin{equation}\label{eqb-hi-basic}
	\sup_{t>0}\bigl(\|\rho(t)\|_{L^\gamma}^\gamma
	+\|\sqrt\rho u(t)\|_{L^2}^2\bigr)
	+\int_0^\infty\|\nabla u\|_{L^2}^2\,dt\le C.
\end{equation}

Choose
$0\le\psi\le1$ in $C_c^\infty(a-1/2,b+1/2)$ with
$\psi=1$ on $[a,b]$ and $\|\psi'\|_{L^1}\le C$.
Now we take a test function 
\[
\Phi_\varepsilon(x,t)
=\psi(t)\mathcal B\bigl([\rho^\theta]_\varepsilon
-\langle[\rho^\theta]_\varepsilon\rangle\bigr).
\]
for the momentum equation $\eqref{eqb-system}_{2}$, which yields
\begin{equation}\label{eqb-hi-identity}
\begin{aligned}
&\int_0^\infty\psi\int_\Omega\rho^\gamma[\rho^\theta]_\varepsilon\,dxdt
+\int_0^\infty\psi'\int_\Omega\rho u\cdot H\,dxdt\\
={}&-\int_0^\infty\psi'\int_\Omega\rho u\cdot(H_\varepsilon-H)\,dxdt
-\int_0^\infty\psi\int_\Omega\rho u\cdot\partial_tH_\varepsilon\,dxdt\\
&-\int_0^\infty\psi\int_\Omega(\rho u\otimes u):\nabla H_\varepsilon\,dxdt
+\mu\int_0^\infty\psi\int_\Omega\nabla u:\nabla H_\varepsilon\,dxdt\\
&+(\lambda+\mu)\int_0^\infty\psi\int_\Omega\operatorname{div}u
\bigl([\rho^\theta]_\varepsilon
-\langle[\rho^\theta]_\varepsilon\rangle\bigr)\,dxdt\\
&-\int_0^\infty\psi\int_\Omega\rho\nabla F\cdot H_\varepsilon\,dxdt
+\int_0^\infty\psi\langle[\rho^\theta]_\varepsilon\rangle
\int_\Omega\rho^\gamma\,dxdt\\
=:{}&-\int_0^\infty\psi'I_{1,\varepsilon}\,dt
+\int_0^\infty\psi\sum_{i=2}^7 I_{i,\varepsilon}\,dt,
\end{aligned}
\end{equation}
where
\[
H=\mathcal B\bigl(\rho^\theta-\langle\rho^\theta\rangle\bigr),
\qquad H_\varepsilon
=\mathcal B\bigl([\rho^\theta]_\varepsilon
-\langle[\rho^\theta]_\varepsilon\rangle\bigr).
\]

Since $\gamma>3/2$, the choice of $\theta$ gives
\begin{equation}\label{eqb-hi-theta-range}
0<\theta\le\frac16<\frac\gamma3,\qquad
\frac1\gamma+\frac16+\frac{\gamma+6\theta}{6\gamma}
\le\frac23+\frac1{2\gamma}<1.
\end{equation}
We next estimate $I_{i,\varepsilon}$ in turn.
It follows from Lemma~\ref{lem:bogovskii}, \eqref{eqb-hi-basic}, and
\eqref{eqb-hi-theta-range} that
\begin{equation}\label{eqb-Phi-bound}
\|H_\varepsilon\|_{W^{1,\gamma/\theta}}+\|H_\varepsilon\|_\infty
\le C\|[\rho^\theta]_\varepsilon\|_{\gamma/\theta}\le C.
\end{equation}
Combining with Lebesgue's theorem, we obtain
\begin{equation}\label{eqb-Phi-convergence}
H_\varepsilon\longrightarrow H\quad\text{in }L^1(a-1/2,b+1/2;L^\infty).
\end{equation}
Using \eqref{eqb-hi-basic} and \eqref{eqb-Phi-convergence}, we obtain
\[
\left|\int_0^\infty\psi'I_{1,\varepsilon}\,dt\right|
\le m^{1/2}\sup_{t>0}\|\sqrt\rho u\|_2
\|\psi'\|_\infty\|H_\varepsilon-H\|_{L^1(a-1/2,b+1/2;L^\infty)}
\longrightarrow0.
\]

For $I_{2,\varepsilon}$, we can apply Lemma~\ref{lem:mollified-continuity} to get
\begin{equation}\label{eqb-hi-time-derivative}
\partial_tH_\varepsilon=-\mathcal B\operatorname{div}
([\rho^\theta]_\varepsilon u)
+(1-\theta)\mathcal B\bigl([\rho^\theta\operatorname{div}u]_\varepsilon
-\langle[\rho^\theta\operatorname{div}u]_\varepsilon\rangle\bigr)
+\mathcal B(r_\varepsilon-\langle r_\varepsilon\rangle),
\end{equation}
By Lemma~\ref{lem:bogovskii}, we obtain
\[
\begin{aligned}
\|\partial_tH_\varepsilon\|_{\frac{6\gamma}{\gamma+6\theta}}
&\le C\bigl(\|[\rho^\theta]_\varepsilon u\|_{\frac{6\gamma}{\gamma+6\theta}}
+\|[\rho^\theta\operatorname{div}u]_\varepsilon\|_{\frac{2\gamma}{\gamma+2\theta}}
+\|r_\varepsilon\|_{\frac{2\gamma}{\gamma+2\theta}}\bigr)\\
&\le C\|\nabla u\|_2+C\|r_\varepsilon\|_{\frac{2\gamma}{\gamma+2\theta}},
\end{aligned}
\]
It follows from this estimate, H\"older's inequality, and
\eqref{eqb-hi-theta-range} that
\[
\begin{aligned}
\left|\int_0^\infty\psi I_{2,\varepsilon}\,dt\right|
&\le C\int_{a-1/2}^{b+1/2}\|\rho\|_\gamma\|u\|_6
\bigl(\|\nabla u\|_2+\|r_\varepsilon\|_{\frac{2\gamma}{\gamma+2\theta}}\bigr)dt\\
&\le C\int_{a-1/2}^{b+1/2}\|\nabla u\|_2^2dt+o_\varepsilon(1).
\end{aligned}
\]
Moreover, by H\"older's inequality, \eqref{eqb-hi-theta-range}, and \eqref{eqb-Phi-bound}, we have
\[
|I_{3,\varepsilon}|
\le C\|\rho\|_\gamma\|u\|_6^2\|\nabla H_\varepsilon\|_{\gamma/\theta}
\le C\|\nabla u\|_2^2.
\]
For the viscous terms, it follows from \eqref{eqb-hi-theta-range} and
\eqref{eqb-Phi-bound} that
\[
\begin{aligned}
|I_{4,\varepsilon}|+|I_{5,\varepsilon}|
&\le C\|\nabla u\|_2
\bigl(\|\nabla H_\varepsilon\|_2
+\|[\rho^\theta]_\varepsilon\|_2\bigr)\\
&\le C\|\nabla u\|_2\le C(1+\|\nabla u\|_2^2).
\end{aligned}
\]
Finally, by H\"older's inequality and \eqref{eqb-hi-basic}, we get
\[
|I_{6,\varepsilon}|+|I_{7,\varepsilon}|
\le m\|\nabla F\|_\infty\|H_\varepsilon\|_\infty
+C\|[\rho^\theta]_\varepsilon\|_{\gamma/\theta}\|\rho\|_\gamma^\gamma\le C,
\]
and
\[
\left|\int_0^\infty\psi'\int_\Omega\rho u\cdot H\,dxdt\right|
\le m^{1/2}\|\psi'\|_1
\sup_{t>0}\bigl(\|\sqrt\rho u\|_2\|H\|_\infty\bigr)\le C.
\]
Combining these estimates with \eqref{eqb-hi-basic} and \eqref{eqb-hi-identity}, we obtain
\[
\begin{aligned}
\int_0^\infty\psi\int_\Omega\rho^\gamma[\rho^\theta]_\varepsilon\,dxdt
&\le C\left(1+b-a+\int_{a-1/2}^{b+1/2}\|\nabla u\|_2^2dt\right)+o_\varepsilon(1)\\
&\le C(1+b-a)+o_\varepsilon(1).
\end{aligned}
\]
Since $[\rho^\theta]_\varepsilon\ge0$ and
$[\rho^\theta]_\varepsilon\to\rho^\theta$ almost everywhere, we can use Fatou's lemma to have
\[
\int_a^b\int_\Omega\rho^{\gamma+\theta}\,dxdt
\le\liminf_{\varepsilon\to0}
\int_0^\infty\psi\int_\Omega\rho^\gamma[\rho^\theta]_\varepsilon\,dxdt
\le C(1+b-a).
\]
All constants depend only on $\Data$.
\end{proof}

We next introduce the property of Young measure representation and its consequence for
strong convergence.

\begin{lemma}\label{lem:young-measure-convergence}
Let $\Omega\subset\R^3$ be a bounded domain, $Q=(0,1)\times\Omega$,
$\gamma>1$, and $\eta>0$. Suppose that $\rho_n\ge0$ and
\[
\sup_n\|\rho_n\|_{L^{\gamma+\eta}(Q)}\le R_1<\infty.
\]
Then, passing to a subsequence if necessary, there exists a measurable
family of probability measures $\nu_{t,x}$ on $[0,\infty)$ with the
following properties:
\begin{enumerate}[label=\textup{(\roman*)},ref=(\roman*)]
\item\label{item:young-representation}
For every $f\in C([0,\infty))$ satisfying
$|f(z)|\le C_f(1+z^\gamma)$,  we have
\[
f(\rho_n)\rightharpoonup\langle\nu_{t,x},f\rangle
\quad\text{in }L^1(Q).
\]
\item\label{item:young-dirac}
If $\nu_{t,x}=\delta_{\rho(t,x)}$ for almost every $(t,x)\in Q$, then
\[
\rho_n\longrightarrow\rho\quad\text{strongly in }L^\gamma(Q).
\]
\end{enumerate}
\end{lemma}
\begin{proof}
For $L>0$, by Chebyshev's inequality, we have
\[
\sup_n|\{\rho_n\ge L\}|
\le R_1^{\gamma+\eta}L^{-(\gamma+\eta)}\longrightarrow0
\quad\text{as }L\to\infty.
\]
Therefore, by \cite[pp.~208--209]{Ball1989}, there exists a subsequence
of $(\rho_n)$, still denoted by $(\rho_n)$, which generates a Young measure
$(\nu_{t,x})$ in the sense of
Definition~\ref{def:young-measure}. For every measurable $E\subset Q$, it follows from
H\"older's inequality that
\[
\int_E\rho_n^\gamma\,dx\,dt
\le R_1^\gamma|E|^{\eta/(\gamma+\eta)},
\]
Moreover, by the growth condition on $f$, we have
\[
\sup_n\int_Q|f(\rho_n)|^{1+\eta/\gamma}\,dx\,dt
\le 2^{\eta/\gamma}C_f^{1+\eta/\gamma}
\bigl(|Q|+R_1^{\gamma+\eta}\bigr)<\infty.
\]
With $\psi(s)=s^{1+\eta/\gamma}$, it follows from
\cite[Remark~3]{Ball1989} that $\{f(\rho_n)\}$ is sequentially weakly
relatively compact in $L^1(Q)$. Hence,
\ref{item:young-representation} holds.

To prove \ref{item:young-dirac}, suppose that
$\nu_{t,x}=\delta_{\rho(t,x)}$. Applying \ref{item:young-representation}
to $f(z)=z^\gamma$ gives $\rho\in L^\gamma(Q)$.
For $L>0$, let $R_L(z)=\min\{z,L\}$. Applying
\ref{item:young-representation} to $R_L$ and $R_L^2$, we obtain
\[
\begin{aligned}
\|R_L(\rho_n)-R_L(\rho)\|_{L^2(Q)}^2
=&\int_Q R_L(\rho_n)^2
-2\int_Q R_L(\rho_n)R_L(\rho)+\int_Q R_L(\rho)^2\longrightarrow0.
\end{aligned}
\]
For $\delta>0$, by Chebyshev's inequality, we have
\[
\begin{aligned}
|\{|\rho_n-\rho|>\delta\}|
&\le |\{\rho_n>L\}|+|\{\rho>L\}|
+\delta^{-2}\|R_L(\rho_n)-R_L(\rho)\|_{L^2(Q)}^2\\
&\le CL^{-\gamma}
+\delta^{-2}\|R_L(\rho_n)-R_L(\rho)\|_{L^2(Q)}^2.
\end{aligned}
\]
Letting first $n\to\infty$ and then $L\to\infty$, we obtain
$\rho_n\to\rho$ in measure on $Q$. The estimate above for
$\int_E\rho_n^\gamma$ implies uniform integrability of
$|\rho_n-\rho|^\gamma$. It follows from Vitali's theorem that
\[
\|\rho_n-\rho\|_{L^\gamma(Q)}\longrightarrow0.
\]
This proves \ref{item:young-dirac}.
\end{proof}

\begin{remark}\label{rem:young-weak-limit}
Sequential weak relative compactness alone does not identify the weak
limit in Lemma~\ref{lem:young-measure-convergence}~\ref{item:young-representation}.
For the subsequence $(\rho_n)$ fixed in its proof, we use a truncation
argument. By \cite[Theorem, \textup{(iii)}]{Ball1989}, the truncated
functions converge weakly-* in $L^\infty(Q)$.
Testing this convergence against $L^\infty(Q)$ functions and using the
uniform $L^{\gamma+\eta}(Q)$ bound to remove the truncation identifies
$\langle\nu_{t,x},f\rangle$ as the only possible weak $L^1(Q)$ limit.
Together with the sequential weak relative compactness from
\cite[Remark~3]{Ball1989}, this proves \ref{item:young-representation}.
\end{remark}

\section{Proof of Theorem~\ref{thm:main}}\label{sec:small-energy-decay}

To prove Theorem~\ref{thm:main}, we first show that the relative energy
becomes sufficiently small after a time depending only on $\Data$.
The key compactness property is that vanishing dissipation on unit time
intervals implies convergence to equilibrium and vanishing time-integrated
relative energy. Combined with the monotonicity of $w$, this property
will give uniform smallness of the relative energy. To establish this
compactness property, we first introduce the operator used to
construct test functions for the momentum equation.

For compactly supported $f\in L^p(\R^3)$ with $1<p<\infty$, define
\[
\mathcal A f:=\nabla\Delta^{-1}f=\nabla(\Gamma*f),
\qquad \Gamma(x):=-\frac1{4\pi|x|}.
\]
Then $\operatorname{div}\mathcal A f=f$ a.e. in $\R^3$.
We next prove the convergence of $\mathcal A[\chi \zeta(\rho_n)]$ and justify
its use in a test function for the momentum equation.

\begin{lemma}\label{lem:local-test-compactness}
Under the assumptions of Theorem~\ref{thm:main}, we can assume
$Q=(0,1)\times\Omega$ and $(\rho_n,u_n)$ are finite energy weak
solutions of \eqref{eqb-system} on $Q$ in the sense of
Definition~\ref{def:weak-solution}. Suppose that, for some $\eta>0$,
\[
\sup_n\bigl(\|\rho_n\|_{L^\infty(0,1;L^\gamma(\Omega))}
+\|\rho_n\|_{L^{\gamma+\eta}(Q)}\bigr)<\infty,
\qquad \|\nabla u_n\|_{L^2(Q)}\longrightarrow0.
\]
Let $\zeta\in C^\infty([0,\infty))$ be bounded with
$z\zeta'(z)-\zeta(z)$ bounded and $\zeta'$ compactly supported. Fix
$\chi\in C_c^\infty(\Omega)$ and a ball $B\supset\overline\Omega$, we set
\[
v_n=\mathcal A[\chi\zeta(\rho_n)],
\]
where $\chi \zeta(\rho_n)$ is extended by zero to $\R^3$.
Then, passing to a subsequence if necessary,
$\zeta(\rho_n)\overset{*}{\rightharpoonup}\tau$ in $L^\infty(Q)$ and
\begin{equation}\label{eq:local-test-limit}
\begin{gathered}
v_n\longrightarrow\mathcal A[\chi\tau]
\quad\text{in }L^\infty(0,1;L^p(B)),\qquad 1\le p<\infty,\\
\partial_tv_n\longrightarrow0
\quad\text{in }L^2(0,1;L^6(\R^3)).
\end{gathered}
\end{equation}
Moreover, we may take $\psi\varphi v_n$ as a test function in the
momentum equation \eqref{eqb-system}$_2$ for every $\psi\in C_c^\infty(0,1)$ and
$\varphi\in C_c^\infty(\Omega)$.
\end{lemma}
\begin{proof}
Write $\beta(z)=z\zeta'(z)-\zeta(z)$.
By Definition~\ref{def:weak-solution}, we can take $b=\zeta$ in
\eqref{eqb-renormalized} to obtain
\begin{equation}\label{eq:repair-tk-cont}
\partial_t \zeta(\rho_n)+\operatorname{div}(\zeta(\rho_n)u_n)
+\beta(\rho_n)\operatorname{div}u_n=0
\quad\text{in }\mathcal D'(Q).
\end{equation}
Since $\zeta$ is bounded and $\chi$ has compact support, we get
\begin{equation}\label{eq:local-cutoff-bound}
\sup_n\|\chi \zeta(\rho_n)\|_{L^\infty(0,1;L^p(\R^3))}\le C_p,
\qquad 1\le p\le\infty.
\end{equation}
By the Riesz-transform estimates \cite[Chapters~II--III]{Stein1970},
we have, for $1<p<\infty$,
\begin{equation}\label{eq:local-riesz-bounds}
\|\nabla\mathcal A f\|_{L^p(\R^3)}\le C_p\|f\|_{L^p(\R^3)},\qquad
\|\mathcal A\operatorname{div}H\|_{L^p(\R^3)}\le C_p\|H\|_{L^p(\R^3)},
\end{equation}
where $f\in L^p(\R^3)$ and $H\in L^p(\R^3;\R^3)$ have compact support.
For $p=2$, we can apply Sobolev inequality to obtain
\begin{equation}\label{eq:local-potential-bound}
\|\mathcal A f\|_{L^6(\R^3)}
\le C\|\nabla\mathcal A f\|_{L^2(\R^3)}
\le C\|f\|_{L^2(\R^3)}.
\end{equation}
Applying \eqref{eq:local-riesz-bounds} and \eqref{eq:local-potential-bound}
to $f=\chi \zeta(\rho_n(t))$ and using \eqref{eq:local-cutoff-bound}, we obtain
\begin{equation}\label{eq:local-vn-bounds}
\|v_n(t)\|_{L^6(\R^3)}\le C,\qquad
\|\nabla v_n(t)\|_{L^p(\R^3)}\le C_p.
\end{equation}
Moreover, by $|\nabla\Gamma(x-y)|\le C|x-y|^{-2}$, we have
\begin{equation}\label{eq:local-vn-pointwise}
\begin{aligned}
|v_n(t,x)|
&\le C\int_{\R^3}\frac{|\chi(y)\zeta(\rho_n(t,y))|}{|x-y|^2}\,dy\\
&\le C\|\chi \zeta(\rho_n(t))\|_{L^\infty}
\int_{|z|<1}|z|^{-2}\,dz+C\|\chi \zeta(\rho_n(t))\|_{L^1}\le C.
\end{aligned}
\end{equation}
Since $B$ is bounded, it follows from \eqref{eq:local-vn-bounds} and
\eqref{eq:local-vn-pointwise} that
\[
\|v_n(t)\|_{W^{1,p}(B)}
\le C_B\|v_n(t)\|_{L^\infty(\R^3)}
+\|\nabla v_n(t)\|_{L^p(\R^3)}\le C_p.
\]
Consequently, for $1<p<\infty$, we have
\begin{equation}\label{eq:repair-vspace}
\|v_n\|_{L^\infty(0,1;L^6(\R^3))}
+\|v_n\|_{L^\infty((0,1)\times\R^3)}
+\|v_n\|_{L^\infty(0,1;W^{1,p}(B))}\le C_p.
\end{equation}
Multiplying \eqref{eq:repair-tk-cont} by $\chi$ and applying $\mathcal A$,
we obtain
\begin{equation}\label{eq:repair-vtime}
\partial_tv_n=-\mathcal A\operatorname{div}(\chi \zeta(\rho_n)u_n)
+\mathcal A[\nabla\chi\cdot \zeta(\rho_n)u_n
-\chi\beta(\rho_n)\operatorname{div}u_n].
\end{equation}
It follows from \eqref{eq:local-riesz-bounds}, \eqref{eq:local-potential-bound}
and \eqref{eq:repair-vtime}, together with Poincar\'e and Sobolev
inequalities, that
\begin{equation}\label{eq:repair-vtime-small}
\|\partial_tv_n\|_{L^2(0,1;L^6(\R^3))}
\le C\|u_n\|_{L^2(0,1;H^1_0(\Omega))}\longrightarrow0.
\end{equation}

Let $m_n=\int_0^1v_n(t)\,dt$. The fundamental theorem of calculus gives
\[
\sup_{t\in[0,1]}\|v_n(t)-m_n\|_{L^6(\R^3)}
\le\|\partial_tv_n\|_{L^2(0,1;L^6(\R^3))}\longrightarrow0.
\]
For $p>3$, it follows from \eqref{eq:repair-vspace} that
the sequence $(m_n)$ is bounded in $W^{1,p}(B)$.
By compact embedding, passing to a subsequence if necessary, we have
\[
m_n\longrightarrow v\quad\text{in }L^6(B),\qquad
v_n\longrightarrow v\quad\text{in }C([0,1];L^6(B)).
\]
By \eqref{eq:repair-vspace} and interpolation, we obtain
\[
v_n\longrightarrow v\quad\text{strongly in }L^\infty(0,1;L^p(B)),
\qquad 1\le p<\infty.
\]
Since $\zeta$ is bounded, passing to a subsequence of $(\rho_n)$
if necessary, we have
\begin{equation}\label{eq:local-density-weak-limit}
\zeta(\rho_n)\overset{*}{\rightharpoonup}\tau
\quad\text{in }L^\infty(Q),\qquad\text{and}\qquad
\chi\zeta(\rho_n)\rightharpoonup\chi\tau
\quad\text{in }L^2(0,1;L^2(\R^3)),
\end{equation}
where $\chi\tau$ is extended by zero outside $\Omega$.
Using the boundedness in \eqref{eq:local-potential-bound}
of $\mathcal A$, we can pass to the weak limit in
\eqref{eq:local-density-weak-limit} to obtain
\[
v_n=\mathcal A[\chi\zeta(\rho_n)]
\rightharpoonup\mathcal A[\chi\tau]
\quad\text{in }L^2(0,1;L^6(\R^3)).
\]
Comparing this with the strong convergence of $(v_n)$ in
$L^2(0,1;L^6(B))$, we conclude by uniqueness of the weak limit that
\[
v=\mathcal A[\chi\tau]\quad\text{a.e. on }(0,1)\times B.
\]
This proves \eqref{eq:local-test-limit}.

We next verify that the weak formulation of the momentum equation
\eqref{eqb-system}$_2$ holds with $\psi\varphi v_n$ as the test function. Choose
\[
p>\max\left\{3,\frac{3\gamma}{2\gamma-3},
\frac{\gamma+\eta}{\eta}\right\}.
\]
By H\"older inequality, we have
\[
\begin{gathered}
\rho_nu_n\in L^2(0,1;L^{6/5}(\Omega)),\qquad
\rho_nu_n\otimes u_n\in L^1(0,1;L^{\frac{3\gamma}{\gamma+3}}(\Omega)),\\
\rho_n^\gamma\in L^{\frac{\gamma+\eta}{\gamma}}(Q),\qquad \nabla u_n\in L^2(Q).
\end{gathered}
\]
For fixed $n$, let $\Phi_{n,\delta}\in C_c^\infty(Q;\R^3)$ be
space-time mollifications of the zero extension of
$\Phi_n=\psi\varphi v_n$ with $\delta>0$ small enough.
By \eqref{eq:repair-vspace} and \eqref{eq:repair-vtime-small}, passing to a
subsequence if necessary, we have that, as $\delta\to0$,
\begin{equation}\label{eq:local-test-approximation}
\begin{gathered}
\partial_t\Phi_{n,\delta}\longrightarrow\partial_t\Phi_n
\quad\text{in }L^2(0,1;L^6(\Omega)),\\
\nabla\Phi_{n,\delta}\longrightarrow\nabla\Phi_n
\quad\text{in }L^p(Q),\qquad
\Phi_{n,\delta}\longrightarrow\Phi_n\quad\text{in }L^{\gamma'}(Q),\\
\Phi_{n,\delta}(t)\longrightarrow\Phi_n(t)
\quad\text{in }W^{1,p}(\Omega)\text{ for a.e. }t,
\end{gathered}
\end{equation}
and
\begin{equation}\label{eq:local-test-uniform-bound}
\sup_\delta\|\Phi_{n,\delta}\|_{L^\infty(0,1;W^{1,p}(\Omega))}<\infty.
\end{equation}
By H\"older's inequality and \eqref{eq:local-test-approximation}, we have
\begin{equation}\label{eq:local-test-time-limit}
\begin{aligned}
&\left|\iint_Q\rho_nu_n\cdot
\partial_t(\Phi_{n,\delta}-\Phi_n)\,dx\,dt\right|
\le\|\rho_nu_n\|_{L^2(0,1;L^{6/5})}
\|\partial_t(\Phi_{n,\delta}-\Phi_n)\|_{L^2(0,1;L^6)}
\longrightarrow0.
\end{aligned}
\end{equation}
For the pressure and viscosity terms, we similarly obtain
\begin{equation}\label{eq:local-test-stress-limit}
\begin{aligned}
&\iint_Q\bigl(\rho_n^\gamma
|\operatorname{div}(\Phi_{n,\delta}-\Phi_n)|
+|\nabla u_n|\,|\nabla(\Phi_{n,\delta}-\Phi_n)|\bigr)\,dx\,dt\\
\le& C\bigl(\|\rho_n^\gamma\|_{L^{\frac{\gamma+\eta}{\gamma}}(Q)}
+\|\nabla u_n\|_{L^2(Q)}\bigr)
\|\nabla(\Phi_{n,\delta}-\Phi_n)\|_{L^p(Q)}\longrightarrow0.
\end{aligned}
\end{equation}
Since $\nabla F\in L^\infty(\Omega)$, it follows from
\eqref{eq:local-test-approximation} that
\begin{equation}\label{eq:local-test-force-limit}
\begin{aligned}
&\left|\iint_Q\rho_n\nabla F\cdot(\Phi_{n,\delta}-\Phi_n)\,dx\,dt\right|
\l& \|\nabla F\|_{L^\infty(\Omega)}\|\rho_n\|_{L^\gamma(Q)}
\|\Phi_{n,\delta}-\Phi_n\|_{L^{\gamma'}(Q)}\longrightarrow0.
\end{aligned}
\end{equation}
For convection term, it follows from H\"older's inequality,
\eqref{eq:local-test-approximation}, and \eqref{eq:local-test-uniform-bound} that
\[
\begin{aligned}
&\left|\int_\Omega\rho_nu_n\otimes u_n:
\nabla(\Phi_{n,\delta}-\Phi_n)\,dx\right|\\
\le& C\|\rho_nu_n\otimes u_n\|_{L^{\frac{3\gamma}{\gamma+3}}(\Omega)}
\|\nabla(\Phi_{n,\delta}-\Phi_n)\|_{L^p(\Omega)}\\
\le& C\|\rho_nu_n\otimes u_n\|_{L^{\frac{3\gamma}{\gamma+3}}(\Omega)}\in L^1(0,1).
\end{aligned}
\]
By \eqref{eq:local-test-approximation} and Lebesgue
theorem, we obtain
\begin{equation}\label{eq:local-test-convection-limit}
\iint_Q\rho_nu_n\otimes u_n:
\nabla(\Phi_{n,\delta}-\Phi_n)\,dx\,dt\longrightarrow0.
\end{equation}
Taking $\Phi_{n,\delta}$ as a test function in
\eqref{eqb-system}$_2$ and letting $\delta\to0$ in
\eqref{eq:local-test-time-limit}--\eqref{eq:local-test-convection-limit},
we may therefore use $\Phi_n=\psi\varphi v_n$ as a test function.
\end{proof}

We now apply Lemma~\ref{lem:local-test-compactness} to prove strong
convergence of the density to $\rho_s$ on unit time intervals where
the viscous dissipation tends to zero.

\begin{lemma}\label{lem:repair-compact}
Under the hypotheses of Theorem~\ref{thm:main}, let
$Q=(0,1)\times\Omega$ and $(\rho_n,u_n)$ be a sequence of finite
energy weak solutions of \eqref{eqb-system} on $Q$ such that
\[
\int_\Omega\rho_n(t,x)\,dx=m
\quad\text{for a.e. }t\in(0,1),\quad n\ge1.
\]
Suppose that, for some $\eta>0$,
\begin{equation}\label{eq:repair-input}
\sup_n\bigl(\|\rho_n\|_{L^\infty(0,1;L^\gamma(\Omega))}
+\|\rho_n\|_{L^{\gamma+\eta}(Q)}\bigr)\le R_1,
\qquad \|\nabla u_n\|_{L^2(Q)}\longrightarrow0.
\end{equation}
Then the sequence $(\rho_n,u_n)$ satisfies
\begin{equation}\label{eq:repair-strong}
\rho_n\longrightarrow \rho_s\quad\text{ in }L^\gamma(Q),
\qquad
\int_0^1 w_n(t)\dd t\longrightarrow0,
\end{equation}
where
$w_n(t)=\int_\Omega(\frac12\rho_n|u_n|^2+G(\rho_n,\rho_s))\dd x$.
\end{lemma}

\begin{proof}
It follows from Poincar\'e and Sobolev inequalities that
\begin{equation}\label{eq:repair-velocity-small}
u_n\longrightarrow0\quad\text{in }L^2(0,1;H^1_0(\Omega))
\cap L^2(0,1;L^6(\Omega)).
\end{equation}
By H\"older's inequality and the Poincar\'e--Sobolev inequalities, we obtain
\begin{equation}\label{eq:repair-product-small}
\begin{aligned}
\|\rho_nu_n\|_{L^2(0,1;L^{\frac{6\gamma}{\gamma+6}}(\Omega))}
&\le C R_1\|u_n\|_{L^2(0,1;L^6)}\longrightarrow0,\\
\|\rho_nu_n\otimes u_n\|_{L^1(0,1;L^{\frac{3\gamma}{\gamma+3}}(\Omega))}
&\le C R_1\|u_n\|_{L^2(0,1;L^6)}^2\longrightarrow0.
\end{aligned}
\end{equation}
In particular, we have
\[
\rho_nu_n\longrightarrow0
\quad\text{in }L^2(0,1;L^{6/5}(\Omega)).
\]
It follows from \eqref{eq:repair-input} and  Banach--Alaoglu theorem that,
passing to a subsequence if necessary, we have
\begin{equation}\label{eq:repair-density-limits}
\rho_n\overset{*}{\rightharpoonup}\bar\rho
\quad\text{in }L^\infty(0,1;L^\gamma(\Omega)),
\qquad
\rho_n^\gamma\rightharpoonup\widetilde\rho
\quad\text{in }L^{\frac{\gamma+\eta}{\gamma}}(Q).
\end{equation}
Passing to the limit in \eqref{eqb-system} against smooth compactly
supported test functions, we obtain
\begin{equation}\label{eq:repair-static-unclosed}
\partial_t\bar\rho=0,\qquad
\nabla\widetilde\rho=\bar\rho\nabla F
\quad\text{in }\mathcal D'(Q).
\end{equation}

Choose $\zeta\in C^\infty([0,\infty))$ such that $\zeta(z)=z$ for $0\le z\le1$,
$0\le \zeta'\le1$, and $\zeta$ is constant for $z\ge2$. Moreover, we further define 
\[
\zeta_k(z)=k\zeta(z/k),\qquad k\in\mathbb N.
\]
By a diagonal argument, passing to a subsequence if necessary, we have
\[
\zeta_k(\rho_n)\overset{*}{\rightharpoonup}\tau_k\quad\text{in }L^\infty(Q),\qquad
\rho_n^\gamma \zeta_k(\rho_n)\rightharpoonup M_k\quad\text{in }L^{\frac{\gamma+\eta}{\gamma}}(Q).
\]
Choose a family of open sets $\{\Omega_j\}_{j\ge1}$ and functions
$\chi_j\in C_c^\infty(\Omega)$ such that
\[
\Omega_j\subset\subset\Omega_{j+1}\subset\subset\Omega,\qquad
\bigcup_{j\ge1}\Omega_j=\Omega,\qquad
\chi_j\equiv1\quad\text{on }\Omega_{j+1}.
\]
For $v_{j,k,n}=\mathcal A[\chi_j\zeta_k(\rho_n)]$,
by Lemma~\ref{lem:local-test-compactness} and a diagonal argument,
passing to a subsequence if necessary, we obtain
\begin{equation}\label{eq:repair-vstrong}
\begin{gathered}
v_{j,k,n}\longrightarrow v_{j,k}=\mathcal A[\chi_j\tau_k]
\quad\text{in }L^\infty(0,1;L^p(B)),\qquad 1\le p<\infty,\\
\partial_tv_{j,k,n}\longrightarrow0
\quad\text{in }L^2(0,1;L^6(\R^3)),\qquad
\operatorname{div}v_{j,k}=\chi_j\tau_k\quad\text{a.e. in }Q,
\end{gathered}
\end{equation}
where $B$ is a fixed ball containing $\overline\Omega$.
Fix $\psi\in C_c^\infty(0,1)$ and $\varphi\in C_c^\infty(\Omega_j)$.
Applying Lemma~\ref{lem:local-test-compactness} again, we may take
$\Phi_n=\psi\varphi v_{j,k,n}$ as a test function in \eqref{eqb-system}$_2$, which implies
\begin{equation}\label{eq:repair-momentum-test}
\begin{aligned}
&\iint_Q\psi\varphi\chi_j\rho_n^\gamma\zeta_k(\rho_n)\\
=&-\iint_Q\Bigl[\psi\rho_n^\gamma\nabla\varphi\cdot v_{j,k,n}
+\rho_nu_n\cdot(\psi'\varphi v_{j,k,n}+\psi\varphi\partial_tv_{j,k,n})\Bigr]\\
&+\iint_Q\psi\Bigl[(\mu\nabla u_n-\rho_nu_n\otimes u_n):\nabla(\varphi v_{j,k,n})
+(\lambda+\mu)\operatorname{div}u_n\operatorname{div}(\varphi v_{j,k,n})
-\varphi\rho_n\nabla F\cdot v_{j,k,n}\Bigr].
\end{aligned}
\end{equation}
By H\"older's inequality, Lemma~\ref{lem:local-test-compactness},
\eqref{eq:repair-input}, \eqref{eq:repair-product-small}, and
\eqref{eq:repair-vstrong}, we obtain
\begin{align*}
&\left|\iint_Q\rho_nu_n\cdot
  \bigl(\psi'\varphi v_{j,k,n}+\psi\varphi\partial_tv_{j,k,n}\bigr)\right|
+\left|\iint_Q\psi(\mu\nabla u_n-\rho_nu_n\otimes u_n):
  \nabla(\varphi v_{j,k,n})\right|\\
&+\left|\iint_Q\psi(\lambda+\mu)\operatorname{div}u_n\,
  \operatorname{div}(\varphi v_{j,k,n})\right|\\
\le& C_{j,k}\Bigl[
\|\rho_nu_n\|_{L^2(0,1;L^{6/5}(\Omega))}
\bigl(1+\|\partial_tv_{j,k,n}\|_{L^2(0,1;L^6(\mathbb R^3))}\bigr)\\
&\quad+\|\rho_nu_n\otimes u_n\|_{L^1(0,1;L^{\frac{3\gamma}{\gamma+3}}(\Omega))}
+\|\nabla u_n\|_{L^2(Q)}\Bigr]\longrightarrow0.
\end{align*}
It follows from \eqref{eq:repair-input}, \eqref{eq:repair-vstrong} and H\"older's inequality that
\begin{align*}
\left|\iint_Q\psi\rho_n^\gamma\nabla\varphi\cdot(v_{j,k,n}-v_{j,k})\right|
\le& C\|\rho_n^\gamma\|_{L^{\frac{\gamma+\eta}{\gamma}}(Q)}
\|v_{j,k,n}-v_{j,k}\|_{L^\infty(0,1;L^{\frac{\gamma+\eta}{\eta}}(B))}\longrightarrow0,\\
\left|\iint_Q\psi\varphi\rho_n\nabla F\cdot(v_{j,k,n}-v_{j,k})\right|
\le& C\|\rho_n\|_{L^\infty(0,1;L^\gamma(\Omega))}
\|v_{j,k,n}-v_{j,k}\|_{L^\infty(0,1;L^{\gamma'}(B))}\longrightarrow0.
\end{align*}
Using \eqref{eq:repair-density-limits}, we have
\begin{align*}
\iint_Q\psi\rho_n^\gamma\nabla\varphi\cdot v_{j,k}
&\longrightarrow\iint_Q\psi\widetilde\rho\nabla\varphi\cdot v_{j,k},\\
\iint_Q\psi\varphi\rho_n\nabla F\cdot v_{j,k}
&\longrightarrow\iint_Q\psi\varphi\bar\rho\nabla F\cdot v_{j,k}.
\end{align*}
Passing to the limit in
\eqref{eq:repair-momentum-test}, we obtain
\begin{equation}\label{eq:repair-first-limit}
\iint_Q\psi\varphi\chi_j M_k
=-\iint_Q\psi\widetilde\rho\nabla\varphi\cdot v_{j,k}
-\iint_Q\psi\varphi\bar\rho\nabla F\cdot v_{j,k}.
\end{equation}
Applying \eqref{eq:repair-static-unclosed} to
$\psi\varphi v_{j,k}$ and using
$\operatorname{div}v_{j,k}=\chi_j\tau_k$ from \eqref{eq:repair-vstrong}, we obtain
\begin{align*}
0
&=\left\langle\nabla\widetilde\rho-\bar\rho\nabla F,
\psi\varphi v_{j,k}\right\rangle_{\mathcal D'(Q),\mathcal D(Q)}\\
&=-\iint_Q\psi\widetilde\rho\,
\operatorname{div}(\varphi v_{j,k})
-\iint_Q\psi\varphi\bar\rho\nabla F\cdot v_{j,k}\\
&=-\iint_Q\psi\widetilde\rho\,
\bigl(\nabla\varphi\cdot v_{j,k}+\varphi\chi_j\tau_k\bigr)
-\iint_Q\psi\varphi\bar\rho\nabla F\cdot v_{j,k},
\end{align*}
which implies
\begin{equation}\label{eq:repair-second-limit}
\iint_Q\psi\varphi\chi_j\widetilde\rho\tau_k
=-\iint_Q\psi\widetilde\rho\nabla\varphi\cdot v_{j,k}
-\iint_Q\psi\varphi\bar\rho\nabla F\cdot v_{j,k}.
\end{equation}
Subtracting \eqref{eq:repair-second-limit} from
\eqref{eq:repair-first-limit} and using
$\chi_j=1$ on $\supp\varphi$, we obtain
\[
M_k=\widetilde\rho\tau_k
\quad\text{a.e. in }(0,1)\times\Omega_j.
\]
Since $\Omega_j\uparrow\Omega$, it follows that
\begin{equation}\label{eq:repair-covariance}
M_k=\widetilde\rho\tau_k\quad\text{a.e. in }Q,
\qquad k\in\mathbb N.
\end{equation}

By Lemma~\ref{lem:young-measure-convergence}~\ref{item:young-representation},
passing to a subsequence if necessary, there exists a Young measure
$\nu_{t,x}$ such that
\[
\bar\rho=\langle\nu_{t,x},z\rangle,\qquad
\widetilde\rho=\langle\nu_{t,x},z^\gamma\rangle,\qquad
\tau_k=\langle\nu_{t,x},\zeta_k(z)\rangle,\qquad
M_k=\langle\nu_{t,x},z^\gamma \zeta_k(z)\rangle, 
\]
where $z\in[0,\infty)$.
It follows from \eqref{eq:repair-covariance} that
\begin{align*}
0&=\langle\nu_{t,x},z^\gamma \zeta_k(z)\rangle
-\langle\nu_{t,x},z^\gamma\rangle\langle\nu_{t,x},\zeta_k(z)\rangle\\
&=\frac12\int_{[0,\infty)^2}
(a^\gamma-b^\gamma)(\zeta_k(a)-\zeta_k(b))\,d\nu_{t,x}(a)\,d\nu_{t,x}(b).
\end{align*}
For every $k\in\mathbb N$, the integrand in the preceding identity satisfies
\[
0\le (a^\gamma-b^\gamma)(\zeta_k(a)-\zeta_k(b))
\le C_k(a^\gamma+b^\gamma),
\qquad a,b\in[0,\infty),
\]
which yields, for a.e.
$(t,x)\in Q$ and every $k\in\mathbb N$,
\[
(a^\gamma-b^\gamma)(\zeta_k(a)-\zeta_k(b))=0
\quad\text{for }(\nu_{t,x}\otimes\nu_{t,x})\text{-a.e. }(a,b)\in[0,\infty)^2.
\]
If $a\ne b$ and $k>\max\{a,b\}$, then
\[
(a^\gamma-b^\gamma)(\zeta_k(a)-\zeta_k(b))
=(a^\gamma-b^\gamma)(a-b)>0.
\]
Thus, for a.e. $(t,x)\in Q$,
\[
(\nu_{t,x}\otimes\nu_{t,x})
\bigl(\{(a,b)\in[0,\infty)^2:a\ne b\}\bigr)=0.
\]
Combining this with
the identity $\bar\rho=\langle\nu_{t,x},z\rangle$, we have
\[
\nu_{t,x}=\delta_{\bar\rho(x)}\quad\text{for a.e. }(t,x)\in Q.
\]
Using Lemma~\ref{lem:young-measure-convergence}~\ref{item:young-dirac}, we get
\[
\rho_n\longrightarrow\bar\rho\quad\text{ in }L^\gamma(Q).
\]
Hence $\widetilde\rho=\bar\rho^\gamma$, which follows from
\eqref{eq:repair-static-unclosed} that
\[
\partial_t\bar\rho=0\quad\text{in }\mathcal D'(Q),\qquad
\nabla\bar\rho^\gamma=\bar\rho\nabla F
\quad\text{in }\mathcal D'(\Omega).
\]
and 
\[
\int_\Omega\bar\rho\,dx=\int_Q\bar\rho\,dx\,dt=m.
\]
The uniqueness in Remark~\ref{rem:equilibrium} therefore gives
$\bar\rho=\rho_s$.

By H\"older's inequality and \eqref{eq:repair-input}, we obtain
\[
\int_0^1\int_\Omega\rho_n|u_n|^2\dd x\dd t
\le\|\rho_n\|_{L^\infty(0,1;L^{3/2})}
\|u_n\|_{L^2(0,1;L^6)}^2\longrightarrow0.
\]
Since $\rho_s\in L^\infty(\Omega)$ and
$\rho_n\to\rho_s$ in $L^\gamma(Q)$, we have
\[
\begin{aligned}
\int_QG(\rho_n,\rho_s)\dd x\dd t
=&\frac1{\gamma-1}\int_Q(\rho_n^\gamma-\rho_s^\gamma)\dd x\dd t-\frac\gamma{\gamma-1}
\int_Q\rho_s^{\gamma-1}(\rho_n-\rho_s)\dd x\dd t\longrightarrow0.
\end{aligned}
\]
More precisely, for every subsequence $(n_\ell)_{\ell\ge1}$ there exists a
further subsequence $(n_{\ell_m})_{m\ge1}$ such that
\[
\int_0^1 w_{n_{\ell_m}}(t)\dd t\longrightarrow0,
\]
which deduce
\[
\int_0^1 w_n(t)\dd t\longrightarrow0.
\]
Together with $\rho_n\to\rho_s$ in $L^\gamma(Q)$, this is exactly
\eqref{eq:repair-strong}.
\end{proof}

For convenience, we set
\[
\mathcal S:=\left\{(\rho,u)\,\middle|\,
\begin{array}{l}
(\rho,u)\text{ is a finite-energy weak solution of \eqref{eqb-system}},\\
\displaystyle\int_\Omega\rho(t,x)\,dx=m\ \text{for a.e. }t\ge0
\end{array}\right\}.
\]
We next prove uniform convergence of the energy $w(t)$.
The proof combines Lemma~\ref{lem:repair-compact} with the energy
inequality and the monotonicity of $w$.
\begin{lemma}\label{thm:repair-uniform}
Under the hypotheses of Theorem~\ref{thm:main}, for every $\varepsilon>0$ there is $T_\varepsilon<\infty$,
depending only on $\Data$ and $\varepsilon$, such that
\[
w(t)\le\varepsilon
\quad\text{for every }t\ge T_\varepsilon,
\qquad\forall(\rho,u)\in\mathcal S.
\]
\end{lemma}
\begin{proof}
Define
\begin{equation}\label{eq:uniform-cost}
\begin{aligned}
\eta_\varepsilon:=\min\biggl\{1,\inf\biggl\{&
\int_b^{b+1}\!\int_\Omega
\bigl[\mu|\nabla u|^2+(\lambda+\mu)(\operatorname{div}u)^2\bigr]dx\,dt
\Big|\,(\rho,u)\in\mathcal S,\ b\ge1,\ w(b+1)>\varepsilon
\biggr\}\biggr\},
\end{aligned}
\end{equation}
and use the convention $\inf\varnothing:=+\infty.$

We claim that $\eta_\varepsilon>0$. Otherwise, there exist
$(\rho_n,u_n)\in\mathcal S$ and $b_n\ge1$ such that
\begin{equation}
\label{eq:uniform-cost-contradiction}
w_n(b_n+1)>\varepsilon,
\qquad
\int_{b_n}^{b_n+1}\!\int_\Omega
\bigl[\mu|\nabla u_n|^2+(\lambda+\mu)(\operatorname{div}u_n)^2\bigr]dx\,dt
\longrightarrow0.
\end{equation}
After the time translation $t\mapsto t+b_n$, we relabel the translated
sequence by $(\rho_n,u_n)$. Thus
\begin{equation}
\label{eq:uniform-translated-smallness}
w_n(1)>\varepsilon,
\qquad
\int_Q\bigl[\mu|\nabla u_n|^2+(\lambda+\mu)(\operatorname{div}u_n)^2\bigr]dx\,dt
\longrightarrow0.
\end{equation}
It follows from Lemma~\ref{eqb-hi} and \eqref{eqb-hi-basic} that
\begin{equation}
\label{eq:uniform-compactness-bounds}
\sup_n\|\rho_n\|_{L^{\gamma+\theta}(Q)}
+\sup_n\|\rho_n\|_{L^\infty(0,1;L^\gamma(\Omega))}<\infty.
\end{equation}
Moreover, \eqref{eq:uniform-translated-smallness} gives
\begin{equation}\label{eq:uniform-gradient-vanishing}
\mu\|\nabla u_n\|_{L^2(Q)}^2
\le
\int_Q\bigl[\mu|\nabla u_n|^2+(\lambda+\mu)(\operatorname{div}u_n)^2\bigr]dx\,dt
\longrightarrow0.
\end{equation}
Consequently, we obtain from Lemma~\ref{lem:repair-compact} and \eqref{eq:uniform-compactness-bounds}, \eqref{eq:uniform-gradient-vanishing} that
\begin{equation}
\label{eq:uniform-energy-vanishing}
\int_0^1w_n(t)\,dt\longrightarrow0.
\end{equation}
On the other hand, by the monotonicity of $w_n$ and
\eqref{eq:uniform-translated-smallness}, one has
\begin{equation}
\int_0^1w_n(t)\,dt\ge w_n(1)>\varepsilon,
\end{equation}
which is a contradiction.

Define
\begin{equation}\label{eq:uniform-time}
T_\varepsilon:=2+\left\lfloor\frac{w_0}{\eta_\varepsilon}\right\rfloor.
\end{equation}
Suppose that $w(T_\varepsilon)>\varepsilon$. By monotonicity of $w$, we have
\[
w(j+1)>\varepsilon,
\qquad j=1,\ldots,T_\varepsilon-1.
\]
It follows from the definition of $\eta_\varepsilon$ that
\[
\int_j^{j+1}\!\int_\Omega
\bigl[\mu|\nabla u|^2+(\lambda+\mu)(\operatorname{div}u)^2\bigr]dx\,dt
\ge\eta_\varepsilon,
\qquad j=1,\ldots,T_\varepsilon-1.
\]
Consequently, using \eqref{eqb-relative-initial-energy} we get
\[
\begin{aligned}
w_0
&\ge\sum_{j=1}^{T_\varepsilon-1}
\int_j^{j+1}\!\int_\Omega
\bigl[\mu|\nabla u|^2+(\lambda+\mu)(\operatorname{div}u)^2\bigr]dx\,dt\\
&\ge (T_\varepsilon-1)\eta_\varepsilon
=\left(1+\left\lfloor\frac{w_0}{\eta_\varepsilon}\right\rfloor\right)\eta_\varepsilon>w_0,
\end{aligned}
\]
which is a contradiction. Thus we have
\[
w(t)\le w(T_\varepsilon)\le\varepsilon,
\qquad t\ge T_\varepsilon,
\]
by monotonicity of $w$.
\end{proof}

Lemma~\ref{thm:repair-uniform} gives a time after which the energy is
small for all weak solutions under consideration. We next use this
smallness to prove exponential decay.
\begin{lemma}\label{lem:small-energy-decay}
There are $\varepsilon_*,c>0$, depending only on $\Data$, such that
for $T\ge1$,
\[
w(t)\le\varepsilon_*\quad\text{for a.e. }t>T
\quad\Longrightarrow\quad
w(t)\le3w_0e^{-c(t-T)}\quad\text{for a.e. }t>T.
\]
\end{lemma}
\begin{proof}
Let $\theta$ be as in Lemma~\ref{eqb-hi} and set
\begin{equation}\label{eq:weighted-field}
\Phi:=\rho_s^{-1}\mathcal B\!\left(
\rho_s^{1-\theta}(\rho^\theta-\rho_s^\theta)
-\left\langle\rho_s^{1-\theta}(\rho^\theta-\rho_s^\theta)\right\rangle
\right).
\end{equation}
By \cite[proof of Proposition~3.1, (3.48)]{KPX}, for every
$0\le\varphi\in C_c^\infty(T,\infty)$, we have
\begin{equation}\label{eq:absorb-decay}
\begin{aligned}
&\int_T^\infty\!\varphi(t)\!\int_\Omega
\rho_s^{-\theta}(\rho^\gamma-\rho_s^\gamma)
(\rho^\theta-\rho_s^\theta)\,dx\,dt
+\int_T^\infty\!\varphi'(t)\!\int_\Omega\rho u\cdot\Phi\,dx\,dt\\
\le&\left[\frac14+\widetilde C\sup_{t>T}
\left(\|\rho^\theta-\rho_s^\theta\|_{L^{\gamma/\theta}(\Omega)}
+\|\rho-\rho_s\|_{L^\gamma(\Omega)}\right)\right]\\
&\times\int_T^\infty\!\varphi(t)\!\int_\Omega
\rho_s^{-\theta}(\rho^\gamma-\rho_s^\gamma)
(\rho^\theta-\rho_s^\theta)\,dx\,dt
+C_0\int_T^\infty\!\varphi(t)\|\nabla u\|_{L^2(\Omega)}^2\,dt,
\end{aligned}
\end{equation}
where $\widetilde C=\widetilde C(\Data)>0$ and $C_0=C_0(\Data)>0$
are independent of $T$ and the particular solution.

With $C_\theta$ and $C_d$ as in Lemma~\ref{eqb-threshold}, we can choose
\[
\varepsilon_*:=\min\left\{1,
\left(\frac{1}{8\widetilde C C_\theta}\right)^{\gamma/\theta},
\left(\frac{1}{8\widetilde C C_d}\right)^{\max\{2,\gamma\}}\right\}>0.
\]
By \eqref{eqb-norm-control} and the hypothesis $w(t)\le\varepsilon_*$,
we have, for a.e. $t>T$,
\[
\begin{aligned}
\widetilde C\left(
\|\rho^\theta-\rho_s^\theta\|_{L^{\gamma/\theta}(\Omega)}
+\|\rho-\rho_s\|_{L^\gamma(\Omega)}\right)
&\le\widetilde C\left(
C_\theta\varepsilon_*^{\theta/\gamma}
+C_d\varepsilon_*^{1/\max\{2,\gamma\}}\right)\le\frac18+\frac18=\frac14.
\end{aligned}
\]
Absorbing the first term on the right-hand side of
\eqref{eq:absorb-decay} into the left-hand side and using
Lemma~\ref{eqb-coercivity}~\ref{item:G-weighted}, we obtain
\begin{equation}\label{eq:potential-dissipation-decay}
\frac{C_a}{2}\int_\Omega G(\rho,\rho_s)\,dx
-\frac{d}{dt}\int_\Omega\rho u\cdot\Phi\,dx
\le C_0\|\nabla u\|_{L^2(\Omega)}^2
\quad\text{in }\mathcal D'(T,\infty).
\end{equation}

To compare the corrected energy with $w$, we use H\"older's inequality,
the Sobolev embedding, Lemma~\ref{lem:bogovskii}~\ref{item:B-bound},
and Lemma~\ref{eqb-coercivity}~\ref{item:G-weighted} to obtain
\[
\begin{aligned}
\left|\int_\Omega\rho u\cdot\Phi\,dx\right|
&\le\|\sqrt\rho u\|_{L^2(\Omega)}
\|\sqrt\rho\|_{L^3(\Omega)}\|\Phi\|_{L^6(\Omega)}\\
&\le C\left(\int_\Omega\rho|u|^2\,dx\right)^{1/2}
\left(\int_\Omega G(\rho,\rho_s)\,dx\right)^{1/2}
\le C_1w(t),
\end{aligned}
\]
where $C_1=C_1(\Data)>0$.
Taking $\delta:=\min\{(2C_1)^{-1},\mu/(2C_0)\}>0$, we have
\begin{equation}\label{eq:corrected-energy-comparison}
\frac12w(t)\le w(t)-\delta\int_\Omega\rho u\cdot\Phi\,dx
\le\frac32w(t).
\end{equation}
Adding $\delta$ times \eqref{eq:potential-dissipation-decay}
to \eqref{eqb-mechanical-distribution}, we obtain
\[
\frac{d}{dt}\left(w-\delta\int_\Omega\rho u\cdot\Phi\,dx\right)
+\frac{\mu}{2}\|\nabla u\|_{L^2(\Omega)}^2
+\frac{\delta C_a}{2}\int_\Omega G(\rho,\rho_s)\,dx
\le0
\quad\text{in }\mathcal D'(T,\infty).
\]
Moreover, by H\"older's inequality, the Sobolev embedding, and
\eqref{eqb-hi-basic}, we have
\[
\frac12\int_\Omega\rho|u|^2\,dx
\le\frac12\|\rho\|_{L^{3/2}(\Omega)}\|u\|_{L^6(\Omega)}^2
\le C_2\|\nabla u\|_{L^2(\Omega)}^2,
\]
where $C_2=C_2(\Data)>0$.
Taking $c:=\frac23\min\{\mu/(2C_2),\delta C_a/2\}>0$ and using
\eqref{eq:corrected-energy-comparison}, we obtain
\[
\frac{d}{dt}\left[
e^{ct}\left(w-\delta\int_\Omega\rho u\cdot\Phi\,dx\right)\right]\le0
\quad\text{in }\mathcal D'(T,\infty).
\]
Thus, for a.e $s,t\in(0,+\infty)$ and $s<t$, it follows from \eqref{eq:corrected-energy-comparison} that
\[
\begin{aligned}
\frac12w(t)
&\le w(t)-\delta\int_\Omega\rho(t)u(t)\cdot\Phi(t)\,dx\\
&\le e^{-c(t-s)}
\left(w(s)-\delta\int_\Omega\rho(s)u(s)\cdot\Phi(s)\,dx\right)\\
&\le\frac32w_0e^{-c(t-s)}.
\end{aligned}
\]
Letting $s\rightarrow T$, we obtain
$w(t)\le3w_0e^{-c(t-T)}$ for a.e. $t>T$.
\end{proof}

We can now combine the uniform energy threshold with the preceding decay lemma.
\begin{proof}[Proof of Theorem~\ref{thm:main}]
If $w_0=0$, it follows from \eqref{eqb-relative-initial-energy} and Lemma~\ref{eqb-coercivity}~\ref{item:G-positive}  that
\[
w(t)=0,\qquad \rho(t)=\rho_s,\qquad u(t)=0
\quad\text{for a.e. }t>0,
\]
so  we just take $T_0=0$ so that the proof is complete. Assume $w_0>0$.
By Lemma~\ref{thm:repair-uniform} with $\varepsilon=\varepsilon_*$ and
$T_0:=T_{\varepsilon_*}$ as in \eqref{eq:uniform-time}, we have
\[
T_0\ge1,\qquad w(t)\le\varepsilon_*\quad(t\ge T_0).
\]
By Lemma~\ref{lem:small-energy-decay}, we obtain \eqref{eq:main-decay}.
By \eqref{eqb-norm-control} and $w(t)\le\varepsilon_*\le1$, we have
\[
\|\rho(t)-\rho_s\|_{L^\gamma(\Omega)}
\le C_dw(t)^\alpha
\le C_d(3w_0)^\alpha e^{-\alpha c(t-T_0)},
\qquad\text{for a.e. }t>T_0,
\]
which proves \eqref{eq:density-rate}.
In particular, by the triangle inequality, we obtain
\[
\|\rho(t)\|_{L^\gamma(\Omega)}
\le\|\rho_s\|_{L^\gamma(\Omega)}
+\|\rho(t)-\rho_s\|_{L^\gamma(\Omega)}
\le\overline{\rho_s}|\Omega|^{1/\gamma}+C_d.
\]
Using this bound, \eqref{eqb-energy}, \eqref{eq:main-decay},
and H\"older's inequality, we have
\[
\begin{aligned}
\|\sqrt\rho u(t)\|_{L^2(\Omega)}
+\|\rho u(t)\|_{L^{2\gamma/(\gamma+1)}(\Omega)}
\le&
\left(1+\|\rho(t)\|_{L^\gamma(\Omega)}^{1/2}\right)
\|\sqrt\rho u(t)\|_{L^2(\Omega)}\\
\le&
\left[1+\left(\overline{\rho_s}|\Omega|^{1/\gamma}+C_d\right)^{1/2}\right]
\sqrt{2w(t)}\\
\le&\sqrt{6w_0}
\left[1+\left(\overline{\rho_s}|\Omega|^{1/\gamma}+C_d\right)^{1/2}\right]
e^{-\frac c2(t-T_0)},
\end{aligned}
\]
which proves \eqref{eq:momentum-rate}.
Finally, by right continuity of $w$, \eqref{eq:main-decay} holds
for every $t>T_0$. Integrating \eqref{eqb-mechanical-distribution}
in time, we obtain
\[
\begin{aligned}
\mu\int_t^\infty\|\nabla u(s)\|_{L^2(\Omega)}^2\,ds
&\le\int_t^\infty\!\int_\Omega
\bigl[\mu|\nabla u|^2+(\lambda+\mu)(\operatorname{div}u)^2\bigr]dx\,ds\\
&\le w(t)\le3w_0e^{-c(t-T_0)},
\end{aligned}
\]
which proves \eqref{eq:tail-dissipation}.
\end{proof}

\section*{Acknowledgements}
This work is supported by the National Natural Science Foundation of China
(No.~12371219), and the Double-Thousand Plan of Jiangxi Province
(No.~jxsq2023201115).

\end{document}